\documentclass[a4]{article}
\usepackage{amsmath}
\usepackage{amssymb}
\usepackage{amsthm}
\usepackage{amscd}
\usepackage{epic}
\newcommand{\ccc}{{{\mathbf C}}}
\newcommand{\nnn}{{{\mathbf N}}}

\newcommand{\rrr}{{{\mathbf R}}}
\newcommand{\zzz}{{{\mathbf Z}}}
\renewcommand{\ggg}{{\frak{g}}}

\newcommand{\hhh}{{\frak{h}}}

\newtheorem{thm}{Theorem}[section]
\newtheorem{prop}{Proposition}[section]
\newtheorem{lemma}{Lemma}[section]
\newtheorem{cor}{Corollary}[section]

\newtheorem{ex}{Example}[section]

\numberwithin{equation}{section}

\begin{document}
\title{Twisted Vertex Operators and $A$-$D$-$E$ Representations}
\author{Minoru Wakimoto \\ \\
{\normalsize Graduate School of Mathematics, Kyushu University}, \\
{\normalsize Fukuoka 812-8581, Japan}}

\maketitle

\begin{abstract}
In this paper, we study the algebra of twisted vertex operators
over an even integral $\zzz_2$-lattice, 
and give a kind of systematic construction of 
fundamental representations for affine Lie algebras of type
$A$, $D$, $E$ with their irreducible decompositions.
\end{abstract}

\section{Introduction}

In this paper, we construct a family of {\it twisted} vertex 
operators associated to an even integral $\zzz_2$-lattice 
which, together with the action of the corresponding Heisenberg
algebra, is invariant under the commutation relations.

For twisted or non-twisted affinization of a finite-dimensional
simple simply-laced Lie algebra $\ggg$, construction of 
fundamental representations is given associated to conjugacy 
classes in the Weyl group of $\ggg$ in \cite{KP1}, 
and associated to automorphisms of $\ggg$ in \cite{Lep1}.
When the $\zzz_2$-lattice is the root $\zzz_2$-lattice of $\ggg$, 
our representation is just the one associated to the automorphism 
of $\ggg$ which is $(-1)$ times the identity transformation on 
the Cartan subalgebra, and so the construction corresponding
to the longest element in the Weyl group in particular 
when $\ggg$ is a simple Lie algebra of type $D_{2m}$, $E_7$ or $E_8$. 
We make a detail analysis on the structure of our representation 
and give its irreducible decomposition explicitly.

This work was motivated by the recent intensive research of 
M. Noumi and Y. Yamada {\it et al} on the Painlev\'e VI equation and 
its Lie algebraic interpretation. 
The author would like to express hearty thanks to 
Professor M. Noumi for private communication and explanation 
on his works at the International Workshop on Integrable Models,
Combinatorics and Representation Theory held in Kyoto on August 2001,
and to Professor E. Date for kind information on $E_8$.

\section{Twisted vertex operators}
\label{section:twistedvo}

Given a positive integer $n$, we consider $\ccc^n$ with a 
non-degenerate symmetric bilinear form $( \,\ | \,\ )$ 
defined by
\begin{equation}
\label{eqn:inner:(1.1)}
(\lambda | \mu) \,\ := \,\ 
\sum^n_{j=1}\lambda_j\mu_j
\end{equation}
for \, $\lambda = (\lambda_1, \cdots, \lambda_n), \, \mu = 
(\mu_1, \cdots, \mu_n) \in \ccc^n$.

Let us consider the space $\ccc [ x^{(j)}_r ; \, j = 1, \cdots, n, \,\ 
r \in \nnn_{\rm odd} ]$. For 
$\lambda = (\lambda_1, \cdots, \lambda_n), \, 
\mu = (\mu_1, \cdots, \mu_n) \in \ccc^n$, we define operators
$U_{\lambda}^{\pm}(z)$, $U_{\lambda}(z)$ and 
$U_{\lambda ; \mu }(z,w)$ on this space as follows:
\begin{subequations}
\begin{equation}
\label{eqn:vertex:(1.1a)}
\begin{array}{cclc}
U_{\lambda}^+(z)
&:=&{\displaystyle
\exp \bigg(
-\sum^{n}_{j=1}\sum_{r \in \nnn_{\rm odd}}\lambda_j
\frac{\partial}{\partial x^{(j)}_r} \frac{z^{-r}}{r}\bigg),} 
& \qquad \quad \cr
U_{\lambda}^-(z)
&:=&{\displaystyle
\exp \bigg(
\sum^{n}_{j=1}\sum_{r \in \nnn_{\rm odd}}\lambda_j x^{(j)}_r z^r\bigg),
} & \qquad \quad \cr
U_{\lambda}(z)&:=&U_{\lambda}^-(z)U_{\lambda}^+(z)
& \qquad \quad 
\end{array}
\end{equation}
and
\begin{eqnarray}
U_{\lambda ; \mu}(z,w)&:=&
U_{\lambda}^-(z)U_{\mu}^-(w)U_{\lambda}^+(z)U_{\mu}^+(w) \cr
&=&
\exp\bigg(
\sum^{n}_{j=1}\sum_{r \in \nnn_{\rm odd}}
(\lambda_jz^r+\mu_jw^r) x^{(j)}_r \bigg) \cr
& & \times \exp \bigg(
-\sum^{n}_{j=1}\sum_{r \in \nnn_{\rm odd}}
\frac{\lambda_jz^{-r}+\mu_jw^{-r}}{r}
\frac{\partial}{\partial x^{(j)}_r} \bigg).
\label{eqn:vertex:(1.1b)} 
\end{eqnarray}
\end{subequations}

From this definition, it follows that
\begin{equation}
\label{eqn:vertex:(1.2)}
\begin{array}{ccc}
U_{\lambda}(-z)&=&U_{-\lambda}(z), \cr
U_{\lambda ; \mu}(-z, w)&=&U_{-\lambda ; \mu}(z, w), \cr
U_{\lambda ; \mu}(z, -w)&=&U_{\lambda ; -\mu}(z, w).
\end{array}
\end{equation}
and that
\begin{equation}
\label{eqn:vertex:(1.3)}
\begin{array}{ccl}
U_{\lambda ; \lambda }(z,-z) &=& 
1 \,\ := \,\ \text{the identity operator}, \cr
{\displaystyle 
\frac{\partial}{\partial z}U_{\lambda ; \lambda }(z,w)}
&=&
{\displaystyle
\sum^{n}_{j=1}\lambda_j
\bigg(\sum_{r \in \nnn_{\rm odd}}rx^{(j)}_rz^{r-1}\bigg)
U_{\lambda ; \lambda }(z,w)}
\cr
& &
{\displaystyle +
U_{\lambda ; \lambda }(z,w)
\sum^{n}_{j=1}\lambda_j
\bigg(\sum_{r \in \nnn_{\rm odd}}z^{-r-1}
\frac{\partial}{\partial x^{(j)}_r} \bigg). }
\end{array}
\end{equation}
Since
$$
U_{\lambda; \mu}(z,z)=U_{\lambda+\mu}(z)
\qquad \text{and} \qquad 
U_{\lambda; \mu}(-z,z)=U_{-\lambda+\mu}(z),
$$
the constant terms of the Taylor series expansion of 
$U_{\lambda; \mu}(z,w)$  around $z=\pm w$ are given as follows:
\begin{equation}
\label{eqn:vertex:(1.7)}
\begin{array}{cclcl}
U_{\lambda; \mu}(z,w)&=&U_{\lambda+\mu}(w) +O(z-w) 
& \qquad & \text{around} \,\ z=w, 
\cr
U_{\lambda; \mu}(z,w)&=&U_{-\lambda+\mu}(w) +O(z+w) 
& \qquad & \text{around} \,\ z=-w. 
\end{array}
\end{equation}

For $j \in \{1, \cdots , n\}$ and $r \in \nnn_{\rm odd}$,
we put
\begin{subequations}
\begin{equation}
\label{eqn:vertex:(1.4a)}
a^{(j)}_r := \frac{\partial}{\partial x^{(j)}_r}, \qquad 
a^{(j)}_{-r} := rx^{(j)}_r \, ,
\end{equation}
and consider the fields
\begin{equation}
\label{eqn:vertex:(1.4b)}
a^{(j)}(z) := \sum_{r \in \zzz_{\rm odd}}a^{(j)}_rz^{-r-1}
\end{equation}
\end{subequations}
for $j=1, \cdots, n$. Notice that
\begin{equation}
\label{eqn:vertex:(1.5)}
a^{(j)}(-z) \,\ = \,\ a^{(j)}(z) \qquad 
\text{for} \,\ j=1, \cdots, n.
\end{equation}
Then from \eqref{eqn:vertex:(1.3)}, one has
\begin{equation}
\label{eqn:vertex:(1.6)}
\left.
\frac{\partial}{\partial z}
U_{\lambda ; \lambda }(z,w) \right\vert_{z=-w}
\, = \, 
\sum^{n}_{j=1}\lambda_ja^{(j)}(w).
\end{equation}

\begin{lemma}
\label{lemma:vertex:(1.1)}
Let $\lambda  \in \ccc^n$. Then the Taylor series expansions 
of $U_{\lambda ; \pm \lambda}(z, w)$ around $z=\mp w$ are 
given as follows:
\begin{enumerate}
\item[{\rm 1)}] \quad ${\displaystyle
U_{\lambda ; \lambda}(z,w) = 1+(z+w)
\sum_{j=1}^n\lambda_ja^{(j)}(w) +O((z+w)^2),}$
\item[{\rm 2)}] \quad ${\displaystyle
U_{\lambda ; -\lambda}(z,w) = 1+(z-w)
\sum_{j=1}^n\lambda_ja^{(j)}(w) +O((z-w)^2).}$
\end{enumerate}
\end{lemma}

\begin{proof}
To prove 1), we compute the Taylor series expansion of \, 
$U_{\lambda ; \lambda }(z,w)$ \, around  $z=-w$, by 
using \eqref{eqn:vertex:(1.3)} and \eqref{eqn:vertex:(1.6)} as follows:
\begin{eqnarray*}
U_{\lambda ; \lambda }(z,w)
&=&
U_{\lambda ; \lambda }(-w,w)
+
\left.
\frac{\partial}{\partial z}U_{\lambda ; \lambda }(z,w) 
\right\vert_{z=-w}
+O((z+w)^2)
\cr
&=&1+ 
\sum^n_{j=1}\lambda_ja^{(j)}(w)+O((z+w)^2).
\end{eqnarray*}
2) follows from 1) and \eqref{eqn:vertex:(1.2)} and 
\eqref{eqn:vertex:(1.5)}.
\end{proof}

\begin{lemma}
\label{lemma:vertex:(1.2)}
For $\lambda, \mu \in \ccc^n$ satisfing 
$(\lambda | \mu) \in 2\zzz$, 
\begin{enumerate}
\item[{\rm 1)}] \quad ${\displaystyle
U_{\lambda}^+(z)U_{\mu}^-(w)
= \iota_{z,w}\left(
\left(\frac{z-w}{z+w}\right)^{(\lambda|\mu)/2} \right)
U_{\mu}^-(w)U_{\lambda}^+(z),}$
\item[{\rm 2)}] \quad ${\displaystyle
U_{\lambda}(z)U_{\mu}(w)
= \iota_{z,w}\left(
\left(\frac{z-w}{z+w}\right)^{(\lambda|\mu)/2}\right)
U_{\lambda ; \mu}(z, w),}$
\end{enumerate}
where $\iota_{z,w}$ means the expansion into the Taylor series in 
the domain $|z| > |w|$.
\end{lemma}

\begin{proof}
For the proof of this lemma, we first notice the commutation 
relation of operators in one variable $x$:
$$
e^{a \frac{\partial}{\partial x}} \circ e^{bx}
\, = \, 
e^{ab}e^{bx} \circ e^{a \frac{\partial}{\partial x}}
$$
for $a, b \in \ccc$, which is easily seen from \, 
$e^{a \frac{\partial}{\partial x}}f(x)  =  f(x+a)$.
Using this, one has
\begin{eqnarray*}
U_{\lambda}^+(z)U_{\mu}^-(w)
&=&
\bigg(
\prod^n_{j=1}
e^{-\lambda_j \mu_j \sum\limits_{r \in \nnn_{\rm odd}}
\frac{z^{-r}w^r}{r}}
\bigg)
U_{\mu}^-(w)U_{\lambda}^+(z)
\cr
&=&
e^{-
\left( \sum\limits^n_{j=1}
\lambda_j \mu_j\right) \sum\limits_{r \in \nnn_{\rm odd}}
\frac{z^{-r}w^r}{r}}
U_{\mu}^-(w)U_{\lambda}^+(z).
\end{eqnarray*}
Then, since $(\lambda|\mu) \in 2\zzz$ by assumption, one has
\begin{eqnarray*}
e^{-(\lambda|\mu)\sum\limits_{r \in \nnn_{\rm odd}}\frac{z^{-r}w^r}{r}}
&=&
\exp \left\{-(\lambda|\mu)\bigg(\sum^{\infty}_{r=1}\frac{z^{-r}w^r}{r}
-
\sum^{\infty}_{r=1}\frac{z^{-2r}w^{2r}}{2r}\bigg)\right\}
\cr
&=&
\iota_{z,w}
\left( \bigg( \frac{1-\frac{w^2}{z^2}}{(1-\frac{w}{z})^2} 
\bigg)^{-(\lambda|\mu)/2} \right)
\cr
&=&
\iota_{z,w}
\left( \bigg( \frac{z+w}{z-w} 
\bigg)^{-(\lambda|\mu)/2} \right),
\end{eqnarray*}
proving 1).  2) follows from 1).
\end{proof}

Note that the vertex operator $U_{\lambda}(z)$ satisfies
\begin{subequations}
\begin{equation}
\label{eqn:vertex:(1.8a)}
\big[
a^{(j)}_r, \, U_{\lambda}(z) \big]
\, = \,
\lambda_jz^rU_{\lambda}(z),
\end{equation}
from which one deduces
\begin{equation}
\label{eqn:vertex:(1.8b)}
\big[
a^{(j)}(z), \, U_{\lambda}(w) \big]
\, = \,
\lambda_j
\sum_{r \in \zzz_{\rm odd}}z^{-r-1}w^r \cdot 
U_{\lambda}(w).
\end{equation}
\end{subequations}
Using the $\delta$-function defined by
\begin{eqnarray}
\delta(z-w) &:=&
\sum_{r \in \zzz}z^{-r-1}w^r
\, = \, 
\sum_{r \in \zzz_{\ge 0}}z^{-r-1}w^r
+\sum_{r \in \zzz_{<0}}z^{-r-1}w^r
\cr
&=&
\big(\iota_{z,w}-\iota_{w,z}\big)
\left(\frac{1}{z-w}\right),
\label{eqn:delta:(1.1)} 
\end{eqnarray}
the formula \eqref{eqn:vertex:(1.8b)} is rewritten as 
\begin{equation}
\label{eqn:vertex:(1.9)}
\big[
a^{(j)}(z), \, U_{\lambda}(w) \big]
\, = \,
\frac{\lambda_j}{2}
\big\{\delta(z-w)-\delta(z+w)\big\}
U_{\lambda}(w),
\end{equation}
since
\begin{subequations}
\begin{eqnarray}
\lefteqn{
\sum_{r \in \zzz_{\rm odd}}z^{-r-1}w^r
\, = \,
\sum_{r \in \zzz}z^{-r-1}w^r
-\sum_{r \in \zzz}z^{-2r-1}w^{2r} } \cr
\cr
&=&
\big(\iota_{z,w}-\iota_{w,z}\big)
\left(\frac{1}{z-w}-\frac{z}{z^2-w^2}\right)
\,\ = \,\ 
\big(\iota_{z,w}-\iota_{w,z}\big)
\left(\frac{w}{z^2-w^2}\right)
\cr
&=&
\big(\iota_{z,w}-\iota_{w,z}\big)
\left(\frac{1/2}{z-w}-\frac{1/2}{z+w}\right).
\label{eqn:delta:(1.2a)}
\end{eqnarray}
Note also that
\begin{eqnarray}
\sum_{r \in \zzz}(-1)^rz^{-r-1}w^r
&=&
\big(\iota_{z,w}-\iota_{w,z}\big)
\left( \frac{1}{z+w}\right),
\label{eqn:delta:(1.2b)} \\
\sum_{r \in \zzz}rz^{-r-1}w^r
&=&
\big(\iota_{z,w}-\iota_{w,z}\big)
\left( \frac{w}{(z-w)^2}\right),
\label{eqn:delta:(1.2c)} \\
\sum_{r \in \zzz}(-1)^rrz^{-r-1}w^r
&=&
\big(\iota_{z,w}-\iota_{w,z}\big)
\left( \frac{-w}{(z+w)^2}\right).
\label{eqn:delta:(1.2d)}
\end{eqnarray}
\end{subequations}

Let $\hhh_{\rrr}$ be an $n$-dimensional real vector space equipped 
with a positive definite symmetric bilinear form $( \,\ | \,\ )$, and 
$Q$ be an even integral $\zzz$-lattice in $\hhh_{\rrr}$ of rank $n$.
We put
$$\Delta := \{ \alpha \in Q \,\ ; \,\ (\alpha | \alpha) =2 \},$$
and decompose $\Delta$ as \, $\Delta = \Delta_+ \cup \Delta_-$.
Let 
$$\nu \,\ : \,\ Q \times Q \,\ \longrightarrow \,\ \{ \pm 1 \}$$
be an asymmetry function (cf. \cite{Kbook1} \S 7.8 and
\cite{Kbook2} \S 5.5); namely a bi-multiplicative
function satisfying the conditions
\begin{equation}
\label{eqn:asymmetry:1}
\begin{array}{ccl}
\nu(\alpha, \alpha) &=&(-1)^{\frac{1}{2}(\alpha | \alpha)}
\cr
\nu(\alpha, \beta) &=&(-1)^{(\alpha | \beta)}\nu(\beta, \alpha)
\end{array}
\end{equation}
for $\alpha, \beta \in Q$. Let  $\ccc\{Q / 2Q\}$  be the 
associative algebra spanned 
over  $\{ e^{\alpha}\}_{\alpha \in Q/2Q}$  with the 
usual non-twisted commutative multiplication, namely 
\begin{equation}
\label{eqn:asymmetry:2}
e^{\alpha}e^{\beta} \,\ := \,\ 
e^{\alpha+\beta}.
\end{equation}
We put
\begin{equation}
\label{eqn:asymmetry:3}
V \,\ := \,\ 
\ccc \{ Q/2Q \}
\otimes 
\ccc \big[x^{(j)}_r \,\ ; \,\ j=1, \cdots, n, \,\ 
r \in \nnn_{\rm odd} \big] .
\end{equation}

Let $\hhh$ be the complexification of $\hhh_{\rrr}$. We extend 
$( \,\ | \,\ )$ to the symmetric bilinear form on $\hhh$ and
fix an orthonormal basis $\{ S_1, \cdots, S_n \}$ of $\hhh$. 
For $\alpha \in \Delta$, we put
$$\alpha (S) := ((\alpha|S_1), \cdots , (\alpha|S_n)) \in \ccc^n,$$
and consider the operator
\begin{equation}
\label{eqn:gamma:1}
\Gamma_{\alpha}(z) := 
\frac{1}{2}e^{\alpha}\nu( \alpha, \,\ \cdot \,\ ) 
\otimes U_{\sqrt{2}\alpha(S)}(z)
\end{equation}
namely
\begin{equation}
\label{eqn:gamma:2}
\Gamma_{\alpha}(z)(e^{\gamma}\otimes f)
:=
\frac{1}{2}\nu(\alpha, \gamma)e^{\alpha+\gamma} \otimes 
U_{(\sqrt{2}(\alpha|S_1), \cdots, \sqrt{2}(\alpha|S_n))}(z)f 
\end{equation}
for $e^{\gamma} \otimes f \in V$.  Notice that 
\begin{equation}
\label{eqn:gamma:3}
\Gamma_{-\alpha}(z) \,\ = \,\ \Gamma_{\alpha}(-z)
\end{equation}
by \eqref{eqn:vertex:(1.2)}.

\begin{lemma}
\label{lemma:gamma:(1.1)}
\,\ For $\alpha, \beta \in Q$, the following formula holds: 
\begin{eqnarray*}
\Gamma_{\alpha}(z)\Gamma_{\beta}(w)(e^{\gamma} \otimes f )
&=&
\frac{1}{4}\nu(\alpha, \beta)\nu(\alpha+\beta, \gamma)
e^{\alpha+\beta+\gamma} 
\cr
&\otimes& \iota_{z,w}\left(
\frac{z-w}{z+w}\right)^{(\alpha|\beta)}
U_{\sqrt{2}\alpha(S); \sqrt{2}\beta(S)}(z, w)f.
\end{eqnarray*}
\end{lemma}

\begin{proof}
For $e^{\gamma}\otimes f \in V$, we have
\begin{eqnarray*}
\lefteqn{
\Gamma_{\alpha}(z)\Gamma_{\beta}(w)(e^{\gamma} \otimes f )
\,\ = \,\ 
\frac{1}{2}\nu(\beta, \gamma)\Gamma_{\alpha}(z)
\left( e^{\beta+\gamma} \otimes U_{\sqrt{2}\beta(S)}(w)f
\right) }
\cr
&=&
\frac{1}{4}
\nu(\beta, \gamma)\nu(\alpha, \beta+\gamma)
e^{\alpha+\beta+\gamma} \otimes 
U_{\sqrt{2}\alpha(S)}(z)
U_{\sqrt{2}\beta(S)}(w)f.
\end{eqnarray*}
Then using Lemma \ref{lemma:vertex:(1.2)} proves the lemma.
\end{proof}

\begin{thm}
\label{thm:gamma:(1.1)}
For $\alpha, \beta \in Q$, the commutators of vertex operators 
$\Gamma_{\alpha}(z)$ and $\Gamma_{\beta}(w)$ are given by the
following formulas:
\begin{enumerate}
\item[{\rm 1)}] \quad If \,\ $(\alpha|\alpha)=2$, then
$$\big[ \Gamma_{\alpha}(z), \, \Gamma_{\alpha}(w)
\big]
=
\big(\iota_{z,w}-\iota_{w,z}\big)
\left(\frac{w}{z+w}
-\frac{w^2}{(z+w)^2}
-\frac{\sqrt{2}w^2}{z+w}
\sum^n_{j=1}(\alpha|S_j)a^{(j)}(w)
\right).
$$
\item[{\rm 2)}] \quad If \,\ $(\alpha|\beta)=1$, then
$$
\big[ \Gamma_{\alpha}(z), \, \Gamma_{\beta}(w)\big] 
\, = \, 
-\nu(\alpha, \beta)
\big(\iota_{z,w}-\iota_{w,z}\big) \left(\frac{w}{z+w}\right)
\Gamma_{-\alpha+\beta}(w).
$$
\item[{\rm 3)}] \quad If \,\ $(\alpha|\beta)=-1$, then
$$
\big[ \Gamma_{\alpha}(z), \, \Gamma_{\beta}(w)\big] 
\, = \, 
\nu(\alpha, \beta)
\big(\iota_{z,w}-\iota_{w,z}\big) \left(\frac{w}{z-w}\right)
\Gamma_{\alpha+\beta}(w).
$$
\item[{\rm 4)}] \quad If \,\ $(\alpha|\beta)=0$, then
\quad $
\big[ \Gamma_{\alpha}(z), \, \Gamma_{\beta}(w)\big] \, = \, 0.$
\end{enumerate}
\end{thm}

\begin{proof}
1) 
Letting $\beta=\alpha$ in Lemma \ref{lemma:gamma:(1.1)} and
using $\nu(\alpha, \alpha)=-1$, one has
$$
\Gamma_{\alpha}(z)\Gamma_{\alpha}(w)(e^{\gamma} \otimes f )
\, = \, 
-\frac{1}{4}e^{\gamma} \otimes 
\iota_{z,w}\left( \frac{z-w}{z+w} \right)^2 
U_{\sqrt{2}\alpha(S) ; \sqrt{2}\alpha(S)}(z,w)f.
$$
Then, exchanging $z$ and $w$, this gives 
\begin{eqnarray*}
\Gamma_{\alpha}(w)\Gamma_{\alpha}(z)(e^{\gamma} \otimes f )
&=&
-\frac{1}{4}e^{\gamma} \otimes 
\iota_{w,z}\left( \frac{w-z}{z+w} \right)^2 
U_{\sqrt{2}\alpha(S) ; \sqrt{2}\alpha(S)}(w,z)f
\cr
&=&
-\frac{1}{4}e^{\gamma} \otimes 
\iota_{w,z}\left( \frac{z-w}{z+w} \right)^2 
U_{\sqrt{2}\alpha(S) ; \sqrt{2}\alpha(S)}(z,w)f.
\end{eqnarray*}
From these two equations, one has
\begin{eqnarray}
\lefteqn{\big[ \Gamma_{\alpha}(z), \, \Gamma_{\alpha}(w)\big]
(e^{\gamma} \otimes f )} \cr
&=&
-\frac{1}{4}e^{\gamma} \otimes 
\big(\iota_{z,w}-\iota_{w,z}\big)\left( \frac{z-w}{z+w} \right)^2 
U_{\sqrt{2}\alpha(S) ; \sqrt{2}\alpha(S)}(z,w)f.
\label{proof:prop:(1.1.1)}
\end{eqnarray}
Notice that
$$
\left(\frac{z-w}{z+w}\right)^2
=
\left(1-\frac{2w}{z+w}\right)^2
= \, 1-\frac{4w}{z+w}+\frac{4w^2}{(z+w)^2},$$
and so 
$$
\big(\iota_{z,w}-\iota_{w,z}\big)
\bigg( \left(\frac{z-w}{z+w}\right)^2 \bigg)
=
\big(\iota_{z,w}-\iota_{w,z}\big)
\left(-\frac{4w}{z+w}+\frac{4w^2}{(z+w)^2}\right).$$
Then, using Lemma \ref{lemma:vertex:(1.1)}.1), 
the formula \eqref{proof:prop:(1.1.1)} is rewritten as follows:
\begin{eqnarray*}
\lefteqn{\big[ \Gamma_{\alpha}(z), \, \Gamma_{\alpha}(w)\big] } \cr
&=&
\big(\iota_{z,w}-\iota_{w,z}\big)
\left(\frac{w}{z+w}-\frac{w^2}{(z+w)^2}\right)
\cr
& &\times
\left\{
1+(z+w)\sqrt{2}\sum^n_{j=1}(\alpha|S_j)a^{(j)}(w)+O((z+w)^2)
\right\}
\cr
&=&
\big(\iota_{z,w}-\iota_{w,z}\big)
\left(\frac{w}{z+w}-\frac{w^2}{(z+w)^2}
-\frac{\sqrt{2}w^2}{z+w}
\sum^n_{j=1}(\alpha|S_j)a^{(j)}(w) \right),
\end{eqnarray*}
proving 1).

2) Applying Lemma \ref{lemma:gamma:(1.1)} to the case 
$(\alpha|\beta)=1$, one has
\begin{eqnarray*}
\lefteqn{
\Gamma_{\alpha}(z)\Gamma_{\beta}(w)(e^{\gamma} \otimes f )
} \cr
&=&
\frac{1}{4}\nu(\alpha, \beta) \nu(\alpha+\beta, \gamma)
e^{\alpha+\beta+\gamma} 
\otimes \iota_{z,w}\left(
\frac{z-w}{z+w}\right)
U_{\sqrt{2}\alpha(S); \sqrt{2}\beta(S)}(z, w)f
\end{eqnarray*}
and, exchanging $z \leftrightarrow w$ and $\alpha 
\leftrightarrow \beta$, also
\begin{eqnarray*}
\lefteqn{
\Gamma_{\beta}(w)\Gamma_{\alpha}(z)(e^{\gamma} \otimes f )
} \cr
&=&
\frac{1}{4}\nu(\beta, \alpha) \nu(\alpha+\beta, \gamma)
e^{\alpha+\beta+\gamma} 
\otimes \iota_{w,z}\left(
\frac{w-z}{z+w}\right)
U_{\sqrt{2}\beta(S); \sqrt{2}\alpha(S)}(w, z)f
\cr
&=&
\frac{1}{4} \nu(\alpha, \beta)\nu(\alpha+\beta, \gamma)
e^{\alpha+\beta+\gamma} \otimes 
\iota_{w,z}\left(\frac{z-w}{z+w}\right)
U_{\sqrt{2}\alpha(S); \sqrt{2}\beta(S)}(z, w)f,
\end{eqnarray*}
since \, $\nu(\beta, \alpha)=-\nu(\alpha, \beta)$.
Then, from these two equations, one has 
\begin{eqnarray}
& &
[ \Gamma_{\alpha}(z), \, \Gamma_{\beta}(w)](e^{\gamma} \otimes f )
\cr
&=&
\frac{1}{4}\nu(\alpha, \beta)\nu(\alpha+\beta, \gamma)
e^{\alpha+\beta+\gamma} \otimes 
\big(\iota_{z,w}-\iota_{w,z}\big)\left(\frac{z-w}{z+w}\right)
U_{\sqrt{2}\alpha(S); \sqrt{2}\beta(S)}(z, w)f.
\cr
& & \label{proof:prop:(1.1.2)}
\end{eqnarray}
Notice that
$$
\frac{z-w}{z+w} \,\ = \,\ 1-\frac{2w}{z+w}$$
and so 
$$
\big(\iota_{z,w}-\iota_{w,z}\big)
\left(\frac{z-w}{z+w} \right)
=
\big(\iota_{z,w}-\iota_{w,z}\big) \left(\frac{-2w}{z+w}\right).$$
Then, by \eqref{eqn:vertex:(1.7)}, 
the formula \eqref{proof:prop:(1.1.2)} is rewritten as follows:
\begin{eqnarray*}
& &
\big[ \Gamma_{\alpha}(z), \, \Gamma_{\beta}(w)\big] 
(e^{\gamma}\otimes f)
\cr
&=&
\frac{-1}{2}\nu(\alpha, \beta)\nu(\alpha+\beta, \gamma)
e^{\alpha+\beta+\gamma} \otimes 
\big(\iota_{z,w}-\iota_{w,z}\big) \left(\frac{w}{z+w}\right)
U_{\sqrt{2}(-\alpha +\beta)}(w)f
\cr
&=&\frac{-1}{2}\nu(\alpha, \beta)
\big(\iota_{z,w}-\iota_{w,z}\big) \left(\frac{w}{z+w}\right)
\nu(-\alpha+\beta, \gamma)
e^{-\alpha+\beta+\gamma} \otimes 
U_{\sqrt{2}(-\alpha +\beta)}(w)f
\cr
&=&
-\nu(\alpha, \beta)
\big(\iota_{z,w}-\iota_{w,z}\big) \left(\frac{w}{z+w}\right)
\Gamma_{-\alpha+\beta}(w)
(e^{\gamma}\otimes f),
\end{eqnarray*}
proving 2).  

3) follows from 2) and \eqref{eqn:gamma:3} since
\begin{eqnarray*}
\big[ \Gamma_{\alpha}(z), \, \Gamma_{\beta}(w)\big] 
&=&
\big[ \Gamma_{\alpha}(z), \, \Gamma_{-\beta}(-w)\big] 
\cr
&=&-\nu(\alpha, \beta)
\big(\iota_{z,w}-\iota_{w,z}\big) \left(\frac{-w}{z-w}\right)
\Gamma_{-\alpha-\beta}(-w)
\cr
&=&\nu(\alpha, \beta)
\big(\iota_{z,w}-\iota_{w,z}\big) \left(\frac{w}{z-w}\right)
\Gamma_{\alpha+\beta}(w).
\end{eqnarray*}
4) follows from Lemma \ref{lemma:gamma:(1.1)} 
and \eqref{eqn:asymmetry:1}.
\end{proof}

Note that formulas in the above theorem are written, in the 
terminology of usual operator products, as follows: 

\begin{cor}
\label{cor:gamma:(1.1)}
For $\alpha, \beta \in Q$, the operator products of vertex
operators $\Gamma_{\alpha}(z)$ and $\Gamma_{\beta}(w)$ are
given by the following formulas:
\begin{enumerate}
\item[{\rm 1)}] \quad If \,\ $(\alpha|\alpha)=2$, then \quad
$$
\Gamma_{\alpha}(z)\Gamma_{\alpha}(-w)
\sim
\frac{-w}{z-w}-\frac{w^2}{(z-w)^2}
-\frac{\sqrt{2}w^2}{z-w}
\sum^n_{j=1}(\alpha|S_j)a^{(j)}(w).
$$
\item[{\rm 2)}] \quad If \,\ $(\alpha|\beta)=1$, then \quad 
${\displaystyle
\Gamma_{\alpha}(z)\Gamma_{\beta}(-w)
\sim
\frac{\nu(\alpha, \beta)w}{z-w} \cdot 
\Gamma_{\alpha -\beta}(w).
}$
\item[{\rm 3)}] \quad If \,\ $(\alpha|\beta)=-1$, then \quad 
${\displaystyle
\Gamma_{\alpha}(z)\Gamma_{\beta}(w)
\sim
\frac{\nu(\alpha, \beta)w}{z-w} \cdot \Gamma_{\alpha +\beta}(w).
}$
\item[{\rm 4)}] \quad If \,\ $(\alpha|\beta)=0$, then \quad 
${\displaystyle
\Gamma_{\alpha}(z)\Gamma_{\beta}(w)
\sim 0.}$
\end{enumerate}
\end{cor}

\section{Twisted affinization of simply-laced Lie algebras}
\label{section:twisted affine}

In this section, we assume that $\ggg$ is a finite-dimensional 
simple Lie algebra of rank $n$ with a symmetric Cartan matrix.
Fix a cartan subalgebra $\hhh$ of $\ggg$, and let $\Delta$ be
the set of all roots of $\ggg$ with respect to $\hhh$ and
$Q$ be the root lattice.  Let
$( \,\ | \,\ )$ be the invariant bilinear form on $\ggg$
normalized by $(\alpha|\alpha)=2$ for all $\alpha \in \Delta$.
For each root $\alpha$, let $\ggg_{\alpha}$ denote the root
space of $\alpha$. 
It is known (cf. \cite{Kbook1} \S 7.8) that, given an asymmetry 
function $\nu : Q \times Q \rightarrow \{ \pm 1 \}$, 
one can choose root vectors $X_{\alpha} \in \ggg_{\alpha}$
satisfying the condition
\begin{subequations}
\begin{equation}
\label{eqn:asymmetry3a}
[X_{\alpha}, \, X_{\beta}] =
\left\{
\begin{array}{lcl}
\nu(\alpha, \beta)X_{\alpha+\beta} & \qquad &
\text{if} \,\ \alpha+\beta \in \Delta \cr
-H_{\alpha} & \qquad & \text{if} \,\ \alpha+\beta = 0
\end{array}
\right.
\end{equation}
for all $\alpha, \beta \in \Delta$, where $H_{\alpha}$ is the 
element in $\hhh$ corresponding to $\alpha$ under the natural 
identification of $\hhh$ with its dual space $\hhh^{\ast}$ via 
the inner product $( \,\ | \,\ )$.
Notice that this condition means
\begin{equation}
\label{eqn:asymmetry3b}
(X_{\alpha} | X_{-\alpha}) \, = \, -1
\qquad \text{for all} \,\ \alpha \in \Delta.
\end{equation}
\end{subequations}

Let $\sigma$ be the automorphism of $\ggg$ such that
\begin{equation}
\label{eqn:asymmetry4}
\begin{array}{ccccl}
\sigma(H) &=& -H &\qquad & \text{for all} \,\  H \in \hhh , \cr
\sigma(X_{\alpha}) &=& X_{-\alpha} &\qquad & \text{for all} \,\  
\alpha \in \Delta.
\end{array}
\end{equation}
We put 
\begin{eqnarray*}
\ggg_{\bar{0}} &:=& \{ X \in \ggg \,\ ; \,\ \sigma (X)=X \} \cr
\ggg_{\bar{1}} &:=& \{ X \in \ggg \,\ ; \,\ \sigma (X)=-X \},
\end{eqnarray*}
and consider the affine Lie algebra
$$
\widehat{\ggg}(\sigma) 
:= \bigg(
\bigoplus_{j \in \zzz}\ggg_{\bar{j} \, {\rm mod} \, 2} \otimes t^j
\bigg) \oplus \ccc K \oplus \ccc d$$
with the Lie bracket
$$
\begin{array}{ccl}
[X \otimes t^j, \, Y \otimes t^k] &:=&
[X, \, Y] \otimes t^{j+k} + \frac{j}{2}(X|Y)\delta_{j+k,0}K ,\cr
[d, \, X \otimes t^j] &:=& jX \otimes t^j, \cr
[K, \, \widehat{\ggg}(\sigma) ] &:=& \{0\},
\end{array}
$$
for $j,k \in \zzz$ and $X \in \ggg_{\bar{j}}$ and 
$Y \in \ggg_{\bar{k}}$.

For each $\alpha \in \Delta$ and $H \in \hhh$, we define the 
fields
\begin{equation}
\label{eqn:(2.1)}
\begin{array}{ccl}
\widetilde{X}_{\alpha}(z)
&:=&
{\displaystyle
\sum_{j \in \zzz}(X_{\alpha}+(-1)^jX_{-\alpha})_{(j)} z^{-j},}
\cr
H(z) &:=& 
{\displaystyle
\sum_{j \in \zzz_{\rm odd}}H_{(j)} z^{-j-1}, }
\end{array}
\end{equation}
where  $X_{(j)} := X \otimes t^j$ for $X \in \ggg$ and $j \in \zzz$
as usual.  Note that
$$
\widetilde{X}_{-\alpha}(z) = \widetilde{X}_{\alpha}(-z)
\qquad \text{and} \qquad  H(-z) = H(z).
$$
Then

\begin{lemma}
\label{lemma:(2.1)}
Let $\alpha, \beta \in \Delta$ and $H \in \hhh$.  Then
\begin{enumerate}
\item[{\rm 1)}] \quad ${\displaystyle
\big[ H(z), \, \widetilde{X}_{\alpha}(w) \big]
\, = \, 
\frac{\alpha(H)}{2}\widetilde{X}_{\alpha}(w) 
\big\{\delta(z-w)-\delta(z+w)\big\},
}$
\item[{\rm 2)}] \quad ${\displaystyle 
\big[
\widetilde{X}_{\alpha}(z), \, 
\widetilde{X}_{\alpha}(w)
\big] } $
$$ = \big(\iota_{z,w}-\iota_{w,z}\big)
\left(\frac{w}{z+w}-\frac{w^2}{(z+w)^2}\right)K
+\big(\iota_{z,w}-\iota_{w,z}\big)
\left(\frac{-2w^2}{z+w}\right)H_{\alpha}(w).
$$
\item[{\rm 3)}] \quad If \,\ $\alpha \ne \pm \beta$, then 
\begin{eqnarray*}
\big[
\widetilde{X}_{\alpha}(z), \, 
\widetilde{X}_{\beta}(w)
\big]
&=& 
\nu(\alpha, \beta)
\big( \iota_{z,w}-\iota_{w,z}\big)\left(\frac{w}{z-w}\right)
\widetilde{X}_{\alpha+\beta}(w)
\cr
&-&
\nu(\alpha, \beta)
\big( \iota_{z,w}-\iota_{w,z}\big)\left(\frac{w}{z+w}\right)
\widetilde{X}_{-\alpha+\beta}(w).
\end{eqnarray*}
\end{enumerate}
\end{lemma}

\begin{proof}
1) is shown as follows:
\begin{eqnarray*}
\lefteqn{
\big[ H(z), \, \widetilde{X}_{\alpha}(w) \big] }
\cr
&=&
\sum_{\substack{j \in \zzz_{\rm odd} \cr
k \in \zzz}}
\big[
H \otimes t^j, \, 
(X_{\alpha}+(-1)^kX_{-\alpha}) \otimes t^k
\big]z^{-j-1}w^{-k}
\cr
&=& \alpha(H)
\sum_{\substack{j \in \zzz_{\rm odd} \cr
k \in \zzz}}
\big(X_{\alpha} -(-1)^kX_{\alpha} \big) \otimes t^{j+k}
z^{-j-1}w^{-k}
\cr
&=& \alpha(H)
\sum_{\substack{j \in \zzz_{\rm odd} \cr
k \in \zzz}}
\big(X_{\alpha} + (-1)^{j+k}X_{\alpha} \big) 
\otimes t^{j+k}w^{-j-k} \cdot z^{-j-1}w^{j}
\cr
\cr
&=& \alpha(H) \widetilde{X}_{\alpha}(w) 
\sum_{j \in \zzz_{\rm odd}}  z^{-j-1}w^{j}
\cr
&=& \frac{\alpha(H)}{2} \widetilde{X}_{\alpha}(w) 
\big\{\delta(z-w)-\delta(z+w)\big\}
\end{eqnarray*}
by \eqref{eqn:delta:(1.2a)}, proving 1).

For the proof of 2) and 3), we first notice the following:
\begin{eqnarray}
& &
\big[
\widetilde{X}_{\alpha}(z), \, 
\widetilde{X}_{\beta}(w)
\big]
\cr
&=&\sum_{j,k \in \zzz} \big[
(X_{\alpha}+(-1)^jX_{-\alpha})\otimes t^j, \, 
(X_{\beta}+(-1)^kX_{-\beta}) \otimes t^k
\big]z^{-j}w^{-k}
\cr
&=&\sum_{j,k \in \zzz} \big[
X_{\alpha}+(-1)^jX_{-\alpha}, \, 
X_{\beta}+(-1)^kX_{-\beta}\big] \otimes t^{j+k}z^{-j}w^{-k}
\cr
& &+K \sum_{j \in \zzz}\frac{j}{2}
(X_{\alpha}+(-1)^jX_{-\alpha}| X_{\beta}+(-1)^jX_{-\beta})
z^{-j}w^j.
\label{proof:lemma:(2.1)}
\end{eqnarray}
Let us consider the case when $\alpha = \beta$.  Since
\begin{eqnarray*}
\lefteqn{
\big[
X_{\alpha}+(-1)^jX_{-\alpha}, \, 
X_{\alpha}+(-1)^kX_{-\alpha}\big] }
\cr
&=& 
(-1)^k\big[ X_{\alpha}, \, X_{-\alpha} \big] 
+
(-1)^j\big[ X_{-\alpha}, \, X_{\alpha} \big] 
\cr
&=&-(-1)^kH_{\alpha}+(-1)^jH_{\alpha}
\, = \, \left\{
\begin{array}{lcl}
2(-1)^jH_{\alpha}  & & \text{if $j+k$ is odd} \cr
0  & & \text{if $j+k$ is even},
\end{array}
\right.
\end{eqnarray*}
and
$$
\big(
X_{\alpha}+(-1)^jX_{-\alpha} \vert 
X_{\alpha}+(-1)^jX_{-\alpha}\big) 
\,\ = \,\ -2(-1)^j
$$
by \eqref{eqn:asymmetry3b}, the formula 
\eqref{proof:lemma:(2.1)} gives
\begin{eqnarray*}
& &
\big[
\widetilde{X}_{\alpha}(z), \, 
\widetilde{X}_{\alpha}(w)
\big]
\cr
&=&2\sum_{\substack{j,k \in \zzz \cr
j+k= {\rm odd}}}
(-1)^jH_{\alpha}  \otimes t^{j+k}z^{-j}w^{-k}
-K \sum_{j \in \zzz}(-1)^jjz^{-j}w^j
\cr
&=&2\sum_{\substack{j,k \in \zzz \cr
j+k= {\rm odd}}}
H_{\alpha}  \otimes t^{j+k}w^{-j-k}\cdot (-1)^jz^{-j}w^{j}
-K \sum_{j \in \zzz}(-1)^jjz^{-j}w^j
\cr
&=&2H_{\alpha}(w) \sum_{j \in \zzz}(-1)^jz^{-j}w^{j}
-K \sum_{j \in \zzz}(-1)^jjz^{-j}w^j
\cr
&=&2wH_{\alpha}(w)
\big(\iota_{z,w}-\iota_{w,z}\big)\left(\frac{z}{z+w}\right)
-K\big(\iota_{z,w}-\iota_{w,z}\big)\left(\frac{-zw}{(z+w)^2}\right)
\cr
&=&
2wH_{\alpha}(w)\big(\iota_{z,w}-\iota_{w,z}\big)
\left(\frac{-w}{z+w}\right)
-K\big(\iota_{z,w}-\iota_{w,z}\big)
\left(\frac{-w}{z+w}+\frac{w^2}{(z+w)^2}
\right)
\end{eqnarray*}
by \eqref{eqn:delta:(1.2b)} and \eqref{eqn:delta:(1.2d)},
proving 2).

In the case $\alpha \ne \pm \beta$, the second term in 
\eqref{proof:lemma:(2.1)} vanishes, so we have 
\begin{eqnarray*}
& &\big[
\widetilde{X}_{\alpha}(z), \, 
\widetilde{X}_{\beta}(w)
\big]
\cr
&=&\nu(\alpha, \beta)\sum_{j,k \in \zzz}
\big\{
X_{\alpha+\beta}+(-1)^{j+k}X_{-\alpha-\beta}\big\}\otimes t^{j+k}
z^{-j}w^{-k}
\cr
& &+\nu(\alpha, \beta)\sum_{j,k \in \zzz}
\big\{
(-1)^kX_{\alpha-\beta}+(-1)^jX_{-\alpha+\beta} \big\} 
\otimes t^{j+k}z^{-j}w^{-k}
\cr
&=&\nu(\alpha, \beta)\sum_{j,k \in \zzz}
\big\{
X_{\alpha+\beta}+(-1)^{j+k}X_{-\alpha-\beta}\big\}\otimes t^{j+k}
w^{-j-k}\cdot z^{-j}w^j
\cr
& &+\nu(\alpha, \beta)\sum_{j,k \in \zzz}(-1)^j
\big\{
(-1)^{j+k}X_{\alpha-\beta}+X_{-\alpha+\beta} \big\} 
\otimes t^{j+k}w^{-j-k} \cdot z^{-j}w^j
\cr
&=&\nu(\alpha, \beta)\widetilde{X}_{\alpha+\beta}(w)
\sum_{j \in \zzz} z^{-j}w^j
+\nu(\alpha, \beta)\widetilde{X}_{-\alpha+\beta}(w)
\sum_{j \in \zzz} (-1)^jz^{-j}w^j
\cr
&=&\nu(\alpha, \beta)
\widetilde{X}_{\alpha+\beta}(w)
\big( \iota_{z,w}-\iota_{w,z}\big)\left(\frac{z}{z-w}\right)
\cr
& &
+\nu(\alpha, \beta)\widetilde{X}_{-\alpha+\beta}(w)
\big( \iota_{z,w}-\iota_{w,z}\big)\left(\frac{z}{z+w}\right),
\end{eqnarray*}
proving 3).
\end{proof}

Noticing that, for $\alpha, \beta \in \Delta$ such that 
$\alpha \ne \pm \beta$, 
$$
\begin{array}{ccl}
\alpha+\beta \in \Delta &\Longleftrightarrow& 
(\alpha|\beta)=-1 \cr
\alpha-\beta \in \Delta &\Longleftrightarrow& 
(\alpha|\beta)=1 \cr
\alpha \pm \beta \not\in \Delta &\Longleftrightarrow& 
(\alpha|\beta)=0 
\end{array}
$$
and comparing Lemma \ref{lemma:(2.1)} with 
Theorem \ref{thm:gamma:(1.1)}, we obtain

\begin{thm}
\label{thm:(2.1)}
Let $\{S_j \}_{j=1, \cdots, n}$ be an orthonormal basis of
$\hhh$. Then the map $\pi : \widehat{\ggg}(\sigma) \longrightarrow
{\rm End}(V)$ defined by
$$
\left\{
\begin{array}{ccccc}
\widetilde{X}_{\alpha}(z) &\longmapsto& \Gamma_{\alpha}(z) 
& & ({}^{\forall} \alpha \in \Delta)
\cr
H(z) &\longmapsto& {\displaystyle 
\frac{1}{\sqrt{2}}\sum^n_{j=1}(H | S_j)a^{(j)}(z)} & & 
({}^\forall H \in \hhh) \cr
K &\longmapsto & 1 \,\ := \,\ \text{the identity operator} & & \cr
d &\longmapsto & {\displaystyle
-L_0 := -\sum^n_{j=1}\sum_{r \in \nnn_{\rm odd}}
rx^{(j)}_r \frac{\partial}{\partial x^{(j)}_r}} & & 
\end{array}
\right.
$$
is a representation of $\widehat{\ggg}(\sigma)$.
\end{thm}

This representation is not irreducible but a sum of finite 
numbers of fundamental representations.  In the next section,
we study its structure and give its irreducible decomposition.
Since the action of $\widehat{\ggg}(\sigma)$ contains 
all $a^{(j)}_r$'s, all singular vectors belong to the subspace 
$\ccc\{Q/2Q\} \otimes 1$ of $V$, which we simply denote by
$\ccc\{Q/2Q\}$. So, in order to get the irreducible decomposition 
of this representation, one needs only to find out singular 
vectors in the space $\ccc\{Q/2Q\}$.

We note also that the transformation $\sigma$ is the 
longest element in the Weyl group of $\ggg$ if $\ggg$ is of type
$D_n$ ($n$ : even) or $E_7$ or $E_8$. In these cases, 
$\widehat{\ggg}(\sigma)$ is a non-twisted affine Lie algebra and 
the representation $\pi$ is a realization of its fundamental 
representation associated to the longest element in the Weyl
group (cf. \cite{KP1}).  Otherwise, $\widehat{\ggg}(\sigma)$ is 
a twisted affine algebra.

\section{Irreducible decomposition for $A$-$D$-$E$ representations}
\label{section:ADE}

Let $\Pi= \{ \alpha_1, \cdots, \alpha_n\}$ denote the set of
simple roots of a finite-dimensional simple Lie algebra $\ggg$ with 
a symmetric Cartan matrix.  Then an asymmetry function $\nu$ is 
determined by $\nu (\alpha_j, \alpha_k)$ $(1 \le j, k \le n)$ by 
its bi-multiplicative property. Then the Dynkin diagram 
of $\Pi$ with orientation corresponds to $\nu$ as follows:  
$$
\nu(\alpha_j, \alpha_k)
=
\left\{
\begin{array}{cl}
1 & \text{if \quad 
\setlength{\unitlength}{1mm}
\begin{picture}(15,8)
\put(0,0){\circle{3}}
\put(13,0){\circle{3}}
\put(1.5,0){\vector(1,0){10}}
\put(0,4){\makebox(0,0){$\alpha_j$}}
\put(13,4){\makebox(0,0){$\alpha_k$}}
\end{picture}
\quad or $\alpha_j$ is not connected with $\alpha_k$}  
\cr
-1 & \text{if} \,\ j=k \,\ \text{or \quad 
\setlength{\unitlength}{1mm}
\begin{picture}(15,8)
\put(0,0){\circle{3}}
\put(13,0){\circle{3}}
\put(11.5,0){\vector(-1,0){10}}
\put(0,4){\makebox(0,0){$\alpha_j$}}
\put(13,4){\makebox(0,0){$\alpha_k$}}
\end{picture}
. }
\end{array}
\right.
$$

For each $\alpha \in Q$, we define the operator 
$\widehat{X}_{\alpha}$ acting on the space $\ccc\{Q/ 2Q\}$ by
\begin{equation}
\label{eqn:(3.1)}
\widehat{X}_{\alpha}(e^{\gamma}) := 
\frac{\nu(\alpha, \gamma)}{2}e^{\alpha+\gamma}.
\end{equation}
In view of \eqref{eqn:gamma:2} and 
Theorem \ref{thm:(2.1)}, one sees that, when $\alpha$ is a root,
this operator $\widehat{X}_{\alpha}$ 
is just the action of $\widetilde{X}_{\alpha}(z)$ to 
the $\ccc \{Q/2Q\}$-component, or more exactly 
$$
\begin{array}{ccl}
\widehat{X}_{\alpha}
&=& (\widetilde{X}_{\alpha})_{(0)} \,\ \text{on} \,\ \ccc\{Q/2Q\}
\cr
&=& \text{the action of} \,\ X_{\alpha}+X_{-\alpha} 
\,\ \text{on} \,\ \ccc\{Q/2Q\}.
\end{array}
$$

For $c_1, \cdots, c_n \in \{\pm 1\}$, we put
\begin{equation}
\label{eqn:(3.2)}
v(c_1, \cdots, c_n) := 
\prod^n_{j=1}(1+ic_je^{\alpha_j}) \,\ \in \,\ \ccc \{Q/2Q \}.
\end{equation}
Then the collection of these elements 
$\{v(c_1, \cdots, c_n)\}_{c_1, \cdots, c_n \in \{\pm 1\}}$  forms 
a basis of $\ccc\{Q/2Q\}$, and the action of 
$\widehat{X}_{\alpha_j}$  on $\ccc\{Q/2Q\}$  is described in terms
of this basis as follows:

\begin{lemma}
\label{lemma:(3.1)}
Let $1 \le j \le n$ and put
$$\{ k_1, \cdots k_s \} \,\ := \,\ 
\{ 1 \le k \le n \,\ ; \,\ 
k \ne j \,\ \text{and} \,\ \nu(\alpha_j, \alpha_k) = -1 \}.$$
Then
$$
2\widehat{X}_{\alpha_j}v(c_1, \cdots, c_n)
=-ic_j
v(c_1, \cdots, -c_{k_1},  \cdots , -c_{k_s}, \cdots, c_n).
$$
\end{lemma}

\begin{proof}
For the proof of this lemma, we notice that
\begin{enumerate}
\item[(i)] \quad 
$2\widehat{X}_{\alpha_j}(1+ic_je^{\alpha_j})
=-ic_j(1+ic_je^{\alpha_j})$,
\item[(ii)] \quad 
$\widehat{X}_{\alpha_j}((1+ic_ke^{\alpha_k})u)
=
(1+i\nu(\alpha_j, \alpha_k)c_ke^{\alpha_k}) 
\cdot \widehat{X}_{\alpha_j}u$ \\
if $k \ne j$ and $u \in \ccc\{Q/2Q\}$.
\end{enumerate}
Actually (i) holds since
$$2\widehat{X}_{\alpha_j}(1+ic_je^{\alpha_j})
=e^{\alpha_j}+ic_j \nu(\alpha_j, \alpha_j)e^{2\alpha_j}
=e^{\alpha_j}-ic_j
=-ic_j(ic_je^{\alpha_j}+1),
$$
and (ii) is shown as follows:
\begin{eqnarray*}
\widehat{X}_{\alpha_j}((1+ic_ke^{\alpha_k})u)
&=&
\widehat{X}_{\alpha_j}(u+ic_ke^{\alpha_k}u)
\,\ = \,\ 
\widehat{X}_{\alpha_j}u+ic_k
\widehat{X}_{\alpha_j}(e^{\alpha_k}u)
\cr
&=&
\widehat{X}_{\alpha_j}u+ic_k
\nu(\alpha_j, \alpha_k)e^{\alpha_k} \cdot 
\widehat{X}_{\alpha_j}u.
\end{eqnarray*}
Then, by the successive use of (ii), one has
\begin{eqnarray*}
2\widehat{X}_{\alpha_j}\left(
\prod^n_{k=1}(1+ic_ke^{\alpha_k})\right)
&=&
2\widehat{X}_{\alpha_j}\bigg(
\prod_{k \ne j}
(1+ic_ke^{\alpha_k}) \cdot 
(1+ic_je^{\alpha_j})
\bigg) \cr
& & \hspace{-25mm} = \,\ 
\prod_{k \ne j}
(1+i\nu(\alpha_j, \alpha_k)c_ke^{\alpha_k}) \cdot 
2\widehat{X}_{\alpha_j}(1+ic_je^{\alpha_j})
\cr
& & \hspace{-25mm} = \,\ 
\prod_{k \ne j}
(1+i\nu(\alpha_j, \alpha_k)c_ke^{\alpha_k}) \cdot 
(-ic_j)(1+ic_je^{\alpha_j}),
\end{eqnarray*}
proving the lemma.
\end{proof}

We note also that
\begin{equation}
\label{eqn:(3.4)}
2\widehat{X}_{\alpha+ \beta}
=
\nu(\alpha, \beta)
\big(2\widehat{X}_{\alpha}\big)
\big(2\widehat{X}_{\beta}\big)
\qquad
\text{for} \,\ \alpha, \beta \in Q,
\end{equation}
since
$$
\widehat{X}_{\alpha+ \beta}(e^{\gamma})
=
\frac{1}{2}\nu(\alpha+\beta, \gamma)e^{\alpha+\beta+\gamma}
$$
and
\begin{eqnarray*}
\widehat{X}_{\alpha}\left(
\widehat{X}_{\beta}(e^{\gamma})\right)
&=&
\frac{1}{2}\nu(\beta, \gamma)
\widehat{X}_{\alpha}(e^{\beta+\gamma})
\, = \, 
\frac{1}{4}
\nu(\beta, \gamma)
\nu(\alpha, \beta+ \gamma)
e^{\alpha+\beta+\gamma}
\cr
&=&
\frac{1}{4}\nu(\alpha, \beta)
\nu(\alpha+ \beta, \gamma)e^{\alpha+\beta+\gamma}.
\end{eqnarray*}
From this formula and Lemma \ref{lemma:(3.1)}, one obtains the
following:

\begin{lemma}
\label{lemma:(3.2)}
\begin{enumerate}
\item[{\rm 1)}] \,\ Let $1 \le j \le p-1$ and $p \le n-2$  in the 
following diagram
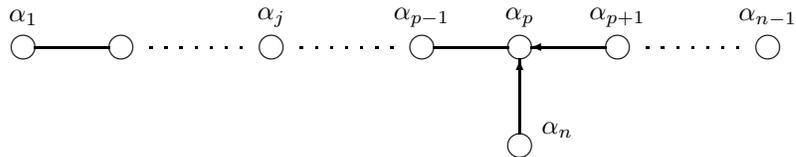
\begin{center}
\setlength{\unitlength}{1mm}
\begin{picture}(103,-20)
\put(1,-4){\circle{3}}
\put(14,-4){\circle{3}}
\put(34,-4){\circle{3}}
\put(54,-4){\circle{3}}
\put(67,-4){\circle{3}}
\put(80,-4){\circle{3}}
\put(100,-4){\circle{3}}
\put(67,-17){\circle{3}}
\put(2.5,-4){\line(1,0){10}}
\put(55.5,-4){\line(1,0){10}}
\put(78.5,-4){\vector(-1,0){10}}
\put(67,-15.5){\vector(0,1){10}}
\put(1,0){\makebox(0,0){$\alpha_1$}}
\put(34,0){\makebox(0,0){$\alpha_j$}}
\put(54,0){\makebox(0,0){$\alpha_{p-1}$}}
\put(67,0){\makebox(0,0){$\alpha_{p}$}}
\put(80,0){\makebox(0,0){$\alpha_{p+1}$}}
\put(100,0){\makebox(0,0){$\alpha_{n-1}$}}
\put(72,-15){\makebox(0,0){$\alpha_{n}$}}
\thicklines
\dottedline{2}(18,-4)(30,-4)
\dottedline{2}(38,-4)(50,-4)
\dottedline{2}(84,-4)(96,-4)
\end{picture}
\end{center}

\vspace{18mm}
\noindent
and \quad 
\setlength{\unitlength}{1mm}
\begin{picture}(15,8)
\put(0,0){\circle{3}}
\put(13,0){\circle{3}}
\put(11.5,0){\vector(-1,0){10}}
\put(0,4){\makebox(0,0){$\alpha_j$}}
\put(13,4){\makebox(0,0){$\alpha_k$}}
\end{picture} .
\,\ Then
\begin{enumerate}
\item[{\rm (i)}] \quad ${\displaystyle
2\widehat{X}_{\alpha_j+\alpha_{p+1}+\alpha_n}
\big(
(1+ic_je^{\alpha_j})
(1+ic_{p+1}e^{\alpha_{p+1}})
(1+ic_ne^{\alpha_n})\big) }$
$$=
ic_jc_{p+1}c_n
(1+ic_je^{\alpha_j})
(1+ic_{p+1}e^{\alpha_{p+1}})
(1+ic_ne^{\alpha_n}),$$
\item[{\rm (ii)}] \quad ${\displaystyle
2\widehat{X}_{\alpha_j+\alpha_{p+1}+\alpha_n}
\big(
(1+ic_je^{\alpha_j})
(1+ic_ke^{\alpha_k})
(1+ic_{p+1}e^{\alpha_{p+1}})
(1+ic_ne^{\alpha_n})\big) }$
$$=
ic_jc_{p+1}c_n
(1+ic_je^{\alpha_j})
(1-ic_ke^{\alpha_k})
(1+ic_{p+1}e^{\alpha_{p+1}})
(1+ic_ne^{\alpha_n}).$$
\end{enumerate}
\item[{\rm 2)}] \,\ Let $1 \le j \le n-4$ in the 
following diagram
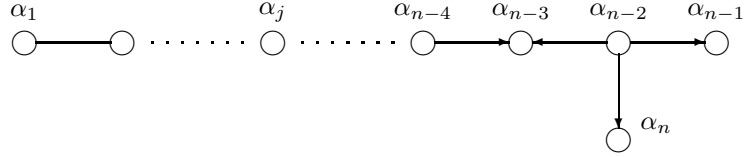
\begin{center}
\setlength{\unitlength}{1mm}
\begin{picture}(100,-20)
\put(1,-4){\circle{3}}
\put(14,-4){\circle{3}}
\put(34,-4){\circle{3}}
\put(54,-4){\circle{3}}
\put(67,-4){\circle{3}}
\put(80,-4){\circle{3}}
\put(93,-4){\circle{3}}
\put(80,-17){\circle{3}}
\put(2.5,-4){\line(1,0){10}}
\put(55.5,-4){\vector(1,0){10}}
\put(78.5,-4){\vector(-1,0){10}}
\put(81.5,-4){\vector(1,0){10}}
\put(80,-5.5){\vector(0,-1){10}}
\put(1,0){\makebox(0,0){$\alpha_1$}}
\put(34,0){\makebox(0,0){$\alpha_j$}}
\put(54,0){\makebox(0,0){$\alpha_{n-4}$}}
\put(67,0){\makebox(0,0){$\alpha_{n-3}$}}
\put(80,0){\makebox(0,0){$\alpha_{n-2}$}}
\put(93,0){\makebox(0,0){$\alpha_{n-1}$}}
\put(85,-15){\makebox(0,0){$\alpha_{n}$}}
\thicklines
\dottedline{2}(18,-4)(30,-4)
\dottedline{2}(38,-4)(50,-4)
\end{picture}
\end{center}

\vspace{18mm}
Then
\begin{enumerate}
\item[{\rm (i)}] \quad ${\displaystyle
2\widehat{X}_{\alpha_j+\alpha_{n-1}+\alpha_n}
\big(
(1+ic_je^{\alpha_j})
(1+ic_{n-2}e^{\alpha_{n-2}})
(1+ic_{n-1}e^{\alpha_{n-1}})
(1+ic_ne^{\alpha_n})\big) }$
$$=
ic_jc_{n-1}c_n
(1+ic_je^{\alpha_j})
(1+ic_{n-2}e^{\alpha_{n-2}})
(1+ic_{n-1}e^{\alpha_{n-1}})
(1+ic_ne^{\alpha_n}),$$
\item[{\rm (ii)}] \quad ${\displaystyle
2\widehat{X}_{\alpha_{n-2}+\alpha_{n-1}+\alpha_n}
\big(
(1+ic_{n-2}e^{\alpha_{n-2}})
(1+ic_{n-1}e^{\alpha_{n-1}})
(1+ic_ne^{\alpha_n})\big) }$
$$=
ic_{n-2}c_{n-1}c_n
(1+ic_{n-2}e^{\alpha_{n-2}})
(1+ic_{n-1}e^{\alpha_{n-1}})
(1+ic_ne^{\alpha_n}).$$
\end{enumerate}
\end{enumerate}
\end{lemma}

We put
\begin{equation}
\label{eqn:(3.3)}
Y_{\alpha} := X_{\alpha}+X_{-\alpha}
\qquad \text{for} \,\ \alpha \in \Delta. 
\end{equation}
Then $Y_{\alpha}$ is an element in $\ggg_{\bar{0}} \cong 
\ggg_{\bar{0}} \otimes t^0 \subset \widehat{\ggg}(\sigma)$ 
and, by the definition \eqref{eqn:(2.1)} of 
field $\widetilde{X}_{\alpha}(z)$, the action of $Y_{\alpha}$'s 
on $\ccc\{Q/2Q\}$ is just equal to $\widehat{X}_{\alpha}$; namely
$$\text{the action of $Y_{\alpha}$ on} \,\ \ccc\{Q/2Q\} 
\,\ = \,\ \widehat{X}_{\alpha}.$$
So one may write $Y_{\alpha}v$  in place 
of $\widehat{X}_{\alpha}v$ for $\alpha \in \Delta$ 
and $v \in \ccc\{Q/2Q\}$.

\subsection{The case $D_n$}

For $D_n$, we consider the following orientation of Dynkin
diagram according as $n$ is even or odd:
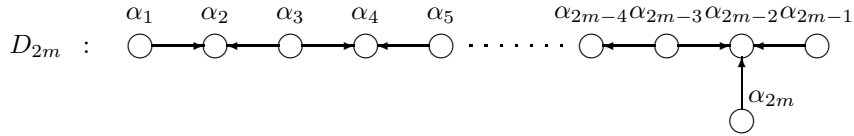
\begin{center}
\setlength{\unitlength}{1mm}
\begin{picture}(110,-18)
\put(15,-4){\circle{3}}
\put(25,-4){\circle{3}}
\put(35,-4){\circle{3}}
\put(45,-4){\circle{3}}
\put(55,-4){\circle{3}}
\put(75,-4){\circle{3}}
\put(85,-4){\circle{3}}
\put(95,-4){\circle{3}}
\put(105,-4){\circle{3}}
\put(95,-14){\circle{3}}
\put(16.5,-4){\vector(1,0){7}}
\put(36.5,-4){\vector(1,0){7}}
\put(86.5,-4){\vector(1,0){7}}
\put(33.5,-4){\vector(-1,0){7}}
\put(53.5,-4){\vector(-1,0){7}}
\put(83.5,-4){\vector(-1,0){7}}
\put(103.5,-4){\vector(-1,0){7}}
\put(95,-12.5){\vector(0,1){7}}
\put(15,0){\makebox(0,0){$\alpha_1$}}
\put(25,0){\makebox(0,0){$\alpha_2$}}
\put(35,0){\makebox(0,0){$\alpha_3$}}
\put(45,0){\makebox(0,0){$\alpha_4$}}
\put(55,0){\makebox(0,0){$\alpha_5$}}
\put(75,0){\makebox(0,0){$\alpha_{2m-4}$}}
\put(85,0){\makebox(0,0){$\alpha_{2m-3}$}}
\put(95,0){\makebox(0,0){$\alpha_{2m-2}$}}
\put(105,0){\makebox(0,0){$\alpha_{2m-1}$}}
\put(99,-11){\makebox(0,0){$\alpha_{2m}$}}
\put(3,-4){\makebox(0,0){$D_{2m} \,\ : $}}
\thicklines
\dottedline{2}(59,-4)(71,-4)
\end{picture}
\end{center}

\vspace{15mm}

\begin{center}
\setlength{\unitlength}{1mm}
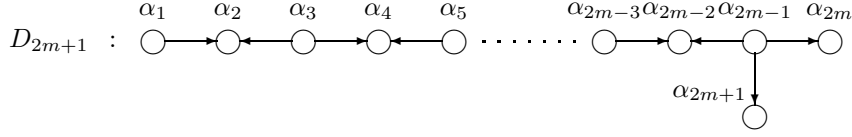
\begin{picture}(110,-18)
\put(15,-4){\circle{3}}
\put(25,-4){\circle{3}}
\put(35,-4){\circle{3}}
\put(45,-4){\circle{3}}
\put(55,-4){\circle{3}}
\put(75,-4){\circle{3}}
\put(85,-4){\circle{3}}
\put(95,-4){\circle{3}}
\put(105,-4){\circle{3}}
\put(95,-14){\circle{3}}
\put(16.5,-4){\vector(1,0){7}}
\put(36.5,-4){\vector(1,0){7}}
\put(76.5,-4){\vector(1,0){7}}
\put(96.5,-4){\vector(1,0){7}}
\put(33.5,-4){\vector(-1,0){7}}
\put(53.5,-4){\vector(-1,0){7}}
\put(93.5,-4){\vector(-1,0){7}}
\put(95,-5.5){\vector(0,-1){7}}
\put(15,0){\makebox(0,0){$\alpha_1$}}
\put(25,0){\makebox(0,0){$\alpha_2$}}
\put(35,0){\makebox(0,0){$\alpha_3$}}
\put(45,0){\makebox(0,0){$\alpha_4$}}
\put(55,0){\makebox(0,0){$\alpha_5$}}
\put(75,0){\makebox(0,0){$\alpha_{2m-3}$}}
\put(85,0){\makebox(0,0){$\alpha_{2m-2}$}}
\put(95,0){\makebox(0,0){$\alpha_{2m-1}$}}
\put(105,0){\makebox(0,0){$\alpha_{2m}$}}
\put(89,-11){\makebox(0,0){$\alpha_{2m+1}$}}
\put(3,-4){\makebox(0,0){$D_{2m+1} \,\ : $}}
\thicklines
\dottedline{2}(59,-4)(71,-4)
\end{picture}
\end{center}

\vspace{18mm}

\begin{prop}
\label{prop:(4.1.1)}
\begin{enumerate}
\item[{\rm 1)}] \quad In the case $n=2m$:
\begin{enumerate}
\item[{\rm (i)}] \quad ${\displaystyle 
2Y_{\alpha_{2j-1}}v(c_1, \cdots, c_{2m})
=-ic_{2j-1}v(c_1, \cdots, c_{2m})}$ \qquad $(1 \le j \le m)$,
\item[{\rm (ii)}] \quad ${\displaystyle 
2Y_{\alpha_{2j}}v(c_1, \cdots, c_{2m}) }$
$$= \left\{
\begin{array}{lcl}
-ic_{2j}v(c_1, \cdots, -c_{2j-1}, c_{2j}, 
-c_{2j+1}, \cdots , c_{2m}) & & (1 \le j \le m-2)
\cr
-ic_{2m-2}v\left(
\begin{array}{rcl}
c_1, \cdots, \, -c_{2m-3}, & \hspace{-2mm} c_{2m-2}, & 
\hspace{-2mm} -c_{2m-1} \cr
& \hspace{-2mm} -c_{2m} & 
\end{array}
\right) & & (j=m-1) \cr
-ic_{2m}v(c_1, \cdots, c_{2m}) & & (j=2m),
\end{array}
\right.
$$
\item[{\rm (iii)}] \quad ${\displaystyle 
2Y_{\alpha_{2j-1}+ 2(\alpha_{2j}+\cdots +
\alpha_{2m-2})+\alpha_{2m-1}+\alpha_{2m}} 
v(c_1, \cdots, c_{2m}) }$
$$= \,\ ic_{2j-1}c_{2m-1}c_{2m}v(c_1, \cdots, c_{2m}) 
\hspace{25mm} (1 \le j \le m-1).$$
\end{enumerate}
\item[{\rm 2)}] \quad In the case $n=2m+1$:
\begin{enumerate}
\item[{\rm (i)}] \quad ${\displaystyle 
2Y_{\alpha_{2j-1}}v(c_1, \cdots, c_{2m+1}) }$
$$
= \left\{
\begin{array}{lcl}
-ic_{2j-1}v(c_1, \cdots, c_{2m+1})
& & (1 \le j \le m)  \cr
-ic_{2m+1}v \left(
\begin{array}{rcl}
c_1, \cdots, c_{2m-2}, & -c_{2m-1}, & c_{2m} \cr
& c_{2m+1} &
\end{array}
\right)
& & (j=m+1),
\end{array}
\right.
$$
\item[{\rm (ii)}] \quad ${\displaystyle 
2Y_{\alpha_{2j}}v(c_1, \cdots, c_{2m+1}) }$
$$= \left\{
\begin{array}{l}
-ic_{2j}v\left(
\begin{array}{rcl}
c_1, \cdots, \, -c_{2j-1}, c_{2j}, \, -c_{2j+1}, \cdots, c_{2m-2}, & 
\hspace{-2mm} c_{2m-1}, & \hspace{-2mm} c_{2m} \cr
& \hspace{-2mm} c_{2m+1} & 
\end{array}
\right) \cr
\hspace{60mm} (1 \le j \le m-1) \cr
-ic_{2m}v\left(
\begin{array}{rcl}
c_1, \cdots, c_{2m-2}, & -c_{2m-1}, & c_{2m} \cr
& c_{2m+1} & 
\end{array}
\right) \qquad (j=m),
\end{array}
\right.
$$
\item[{\rm (iii)}] \quad ${\displaystyle 
2Y_{\alpha_{2j-1}+ 2(\alpha_{2j}+\cdots +
\alpha_{2m-1})+\alpha_{2m}+\alpha_{2m+1}} 
v(c_1, \cdots, c_{2m+1}) }$
$$= \,\ ic_{2j-1}c_{2m}c_{2m+1}v(c_1, \cdots, c_{2m+1}) 
\qquad (1 \le j \le m-1),$$
\item[{\rm (iv)}] \quad ${\displaystyle 
2Y_{\alpha_{2m-1}+\alpha_{2m}+\alpha_{2m+1}} 
v(c_1, \cdots, c_{2m+1}) }$
$$= \,\ ic_{2m-1}c_{2m}c_{2m+1}v(c_1, \cdots, c_{2m+1}) .$$
\end{enumerate}
\end{enumerate}
\end{prop}

For an explicit description of Chevalley generators 
of $\widehat{\ggg}(\sigma)$, 
we define elements $Z_{j,k}$ and $Z_{j,k}^{\prime}$ 
in $\ggg_{\bar{0}}$ for $1 \le j \le k \le n-1$ as follows:
$$
Z_{j,k}
:=
\left\{
\begin{array}{ll}
Y_{\alpha_j+ \cdots + \alpha_k}
+
Y_{\alpha_j+ \cdots + \alpha_k+2(\alpha_{k+1}+\cdots +\alpha_{n-2})
+\alpha_{n-1}+\alpha_{n}} 
& (k \le n-3) \cr
Y_{\alpha_j+ \cdots + \alpha_{n-2}}
+
Y_{\alpha_j+ \cdots + \alpha_{n-2}+\alpha_{n-1}+\alpha_{n}} 
& (k = n-2) \cr
Y_{\alpha_j+ \cdots +\alpha_{n-2}+\alpha_{n-1}}
-
Y_{\alpha_j+ \cdots +\alpha_{n-2}+\alpha_{n}}  & (j<k = n-1) \cr
Y_{\alpha_{n-1}}-Y_{\alpha_{n}} & (j=k = n-1) 
\end{array}
\right.
$$
and
$$
Z_{j,k}^{\prime}
:=
\left\{
\begin{array}{ll}
Y_{\alpha_j+ \cdots + \alpha_k}
-
Y_{\alpha_j+ \cdots + \alpha_k+2(\alpha_{k+1}+\cdots +\alpha_{n-2})
+\alpha_{n-1}+\alpha_{n}}  
& (k \le n-3) \cr
Y_{\alpha_j+ \cdots + \alpha_{n-2}}
-
Y_{\alpha_j+ \cdots + \alpha_{n-2}+\alpha_{n-1}+\alpha_{n}} 
& (k = n-2) \cr
Y_{\alpha_j+ \cdots +\alpha_{n-2}+\alpha_{n-1}}
+
Y_{\alpha_j+ \cdots +\alpha_{n-2}+\alpha_{n}} & (j<k = n-1) \cr
Y_{\alpha_{n-1}}+Y_{\alpha_{n}} & (j=k = n-1). 
\end{array}
\right.
$$
For simplicity, we write  $Z_j := Z_{j,j}$  and  $Z_j^{\prime} 
:= Z_{j,j}^{\prime}$.
Then, by an easy calculation, one can check the following:

\begin{lemma}
\label{lemma:(4.1.1)}
\begin{enumerate}
\item[{\rm 1)}] \quad 
$\big[Z_{j,k}, \, Z_{r,s}^{\prime}\big] = 0 $
\qquad for all $j,k,r,s$.
\item[{\rm 2)}]
\begin{enumerate}
\item[{\rm (i)}] \quad
$\big[Z_{j,r}, \, Z_{j,s}\big] = \left\{
\begin{array}{ccl}
-2\nu(\alpha_r, \alpha_{r+1})Z_{r+1, s} & & (r<s), \cr
2\nu(\alpha_s, \alpha_{s+1})Z_{s+1, r} & & (s<r),
\end{array} \right.$
\item[{\rm (ii)}] \quad
$\big[Z_{j,r}^{\prime}, \, Z_{j,s}^{\prime}\big] = \left\{
\begin{array}{ccl}
-2\nu(\alpha_r, \alpha_{r+1})Z_{r+1, s}^{\prime} & & (r<s), \cr
2\nu(\alpha_s, \alpha_{s+1})Z_{s+1, r}^{\prime} & & (s<r).
\end{array} \right. $
\end{enumerate}
\item[{\rm 3)}]
\begin{enumerate}
\item[{\rm (i)}] \quad
$\big[Z_{j,r}, \, Z_{k,r}\big] = \left\{
\begin{array}{ccl}
-2\nu(\alpha_{k-1}, \alpha_{k})Z_{j, k-1} & & (j<k), \cr
2\nu(\alpha_{j-1}, \alpha_{j})Z_{k, j-1} & & (k<j),
\end{array} \right. $
\item[{\rm (ii)}] \quad
$\big[Z_{j,r}^{\prime}, \, Z_{k,r}^{\prime}\big] = \left\{
\begin{array}{ccl}
-2\nu(\alpha_{k-1}, \alpha_{k})Z_{j, k-1}^{\prime} & & (j<k), \cr
2\nu(\alpha_{j-1}, \alpha_{j})Z_{k, j-1}^{\prime} & & (k<j).
\end{array} \right. $
\end{enumerate}
\item[{\rm 4)}]
\begin{enumerate}
\item[{\rm (i)}] \quad $\big[Z_{j,k-1}, \, Z_{k,s}\big] = 
2\nu(\alpha_{k-1}, \alpha_{k})Z_{j, s},$
\item[{\rm (ii)}] \quad
$\big[Z_{j,k-1}^{\prime}, \, Z_{k,s}^{\prime}\big] = 
2\nu(\alpha_{k-1}, \alpha_{k})Z_{j, s}^{\prime}. $
\end{enumerate}
\item[{\rm 5)}]
\begin{enumerate}
\item[{\rm (i)}] \quad $\big[Z_{j,r}, \, Z_{k,j-1}\big] = 
-2\nu(\alpha_{j-1}, \alpha_{j})Z_{k, r},$
\item[{\rm (ii)}] \quad
$\big[Z_{j,r}^{\prime}, \, Z_{k,j-1}^{\prime}\big] = 
-2\nu(\alpha_{j-1}, \alpha_{j})Z_{k, r}^{\prime}.$
\end{enumerate}
\end{enumerate}
\end{lemma}

For $\ggg=D_n=so(2n)$, the $\sigma$-fixed subalgebra 
$\ggg_{\bar{0}}$ is $so(n) \oplus so(n)$, which is a semisimple 
Lie algebra of type $D_m \oplus D_m$ if $n=2m$ 
and $B_m \oplus B_m$ if $n=2m+1$, where $D_2:=A_1 \oplus A_1$. 
We define elements 
$\widetilde{e}_j, \widetilde{f}_j, \widetilde{h}_j$ 
for $0 \le j \le m$ and 
$\widetilde{e}_j^{\prime}, \widetilde{f}_j^{\prime}, 
\widetilde{h}_j^{\prime}$ for $1 \le j \le m$ 
in $\widehat{\ggg}(\sigma)$ as follows:

\noindent
In the case $n=2m$:
\begin{eqnarray}
& & \left\{
\begin{array}{ccl}
\widetilde{e}_j &:=& \frac{1}{2}
\big\{
Z_{2j-1, 2j}-Z_{2j, 2j+1}-iZ_{2j-1, 2j+1}-iZ_{2j}
\big\} \cr
& & \hspace{50mm} (1 \le j \le m-1) \cr
\widetilde{e}_m &:=& \frac{1}{2}
\big\{
Z_{2m-3, 2m-2}+Z_{2m-2, 2m-1}+iZ_{2m-3, 2m-1}-iZ_{2m-2}
\big\} \cr
\widetilde{e}_j^{\prime} &:=& \frac{1}{2}
\big\{ Z_{2j-1, 2j}^{\prime}-Z_{2j, 2j+1}^{\prime}
-iZ_{2j-1, 2j+1}^{\prime}-iZ_{2j}^{\prime} \big\} \cr
& & \hspace{50mm} (1 \le j \le m-1) \cr
\widetilde{e}_m^{\prime} &:=& \frac{1}{2}
\big\{ Z_{2m-3, 2m-2}^{\prime}+Z_{2m-2, 2m-1}^{\prime}
+iZ_{2m-3, 2m-1}^{\prime}-iZ_{2m-2}^{\prime} \big\} 
\end{array}
\right. \cr
& & \left\{
\begin{array}{ccl}
\widetilde{f}_j &:=& \frac{1}{2}
\big\{
-Z_{2j-1, 2j}+Z_{2j, 2j+1}-iZ_{2j-1, 2j+1}-iZ_{2j} 
\big\} \cr
& & \hspace{50mm} (1 \le j \le m-1) \cr
\widetilde{f}_m &:=& \frac{1}{2}
\big\{
-Z_{2m-3, 2m-2}-Z_{2m-2, 2m-1}+iZ_{2m-3, 2m-1}-iZ_{2m-2}
\big\} \cr
\widetilde{f}_j^{\prime} &:=& \frac{1}{2}
\big\{ -Z_{2j-1, 2j}^{\prime}+Z_{2j, 2j+1}^{\prime}
-iZ_{2j-1, 2j+1}^{\prime}-iZ_{2j}^{\prime} \big\} \cr
& & \hspace{50mm} (1 \le j \le m-1) \cr
\widetilde{f}_m^{\prime} &:=& \frac{1}{2}
\big\{ -Z_{2m-3, 2m-2}^{\prime}-Z_{2m-2, 2m-1}^{\prime}
+iZ_{2m-3, 2m-1}^{\prime}-iZ_{2m-2}^{\prime} \big\} 
\end{array}
\right. \cr
& & \left\{
\begin{array}{cclcl}
\widetilde{h}_j &:=& \frac{i}{2}\big\{ Z_{2j-1}-Z_{2j+1} \big\} 
& \qquad &(1 \le j \le m-1) \cr
\widetilde{h}_m &:=& \frac{i}{2}\big\{ Z_{2m-3}+Z_{2m-1} \big\} 
& & \cr
\widetilde{h}_j^{\prime} &:=& \frac{i}{2}
\big\{ Z_{2j-1}^{\prime}-Z_{2j+1}^{\prime} \big\} 
&\qquad &(1 \le j \le m-1) \cr
\widetilde{h}_m^{\prime} &:=& \frac{i}{2}
\big\{ Z_{2m-3}^{\prime}+Z_{2m-1}^{\prime} \big\} 
& & 
\end{array}
\right. \cr
& & \left\{
\begin{array}{ccl}
\widetilde{e}_0 &:=& \frac{1}{2}
\big\{ i(X_{\alpha_1}-X_{-\alpha_1})+H_{\alpha_1} \big\} 
\otimes t \cr
\widetilde{f}_0 &:=& \frac{1}{2}
\big\{ -i(X_{\alpha_1}-X_{-\alpha_1})+H_{\alpha_1} \big\} 
\otimes t^{-1} \cr
\widetilde{h}_0 &:=& -iY_{\alpha_1}+\frac{K}{2}.
\end{array}
\right.
\label{eqn:(4.1.1)}
\end{eqnarray}
In the case $n=2m+1$:
\begin{eqnarray}
& & \left\{
\begin{array}{ccl}
\widetilde{e}_j &:=& \frac{1}{2}
\big\{
Z_{2j-1, 2j}-Z_{2j, 2j+1}-iZ_{2j-1, 2j+1}-iZ_{2j}
\big\} \cr
& & \hspace{50mm} (1 \le j \le m-1) \cr
\widetilde{e}_m &:=& \frac{1}{2}
\big\{
Z_{2m-2, 2m-1}-iZ_{2m-2}
\big\} \cr
\widetilde{e}_j^{\prime} &:=& \frac{1}{2}
\big\{ Z_{2j-1, 2j}^{\prime}-Z_{2j, 2j+1}^{\prime}
-iZ_{2j-1, 2j+1}^{\prime}-iZ_{2j}^{\prime} \big\} \cr
& & \hspace{50mm} (1 \le j \le m-1) \cr
\widetilde{e}_m^{\prime} &:=& \frac{1}{2}
\big\{ Z_{2m-2, 2m-1}^{\prime}-iZ_{2m-2}^{\prime} \big\} 
\end{array}
\right. \cr
& & \left\{
\begin{array}{ccl}
\widetilde{f}_j &:=& \frac{1}{2}
\big\{
-Z_{2j-1, 2j}+Z_{2j, 2j+1}-iZ_{2j-1, 2j+1}-iZ_{2j}
\big\} \cr
& & \hspace{50mm} (1 \le j \le m-1) \cr
\widetilde{f}_m &:=& \frac{1}{2}
\big\{
-Z_{2m-2, 2m-1}
-iZ_{2m-2}
\big\} \cr
\widetilde{f}_j^{\prime} &:=& \frac{1}{2}
\big\{
-Z_{2j-1, 2j}^{\prime}
+Z_{2j, 2j+1}^{\prime}
-iZ_{2j-1, 2j+1}^{\prime}
-iZ_{2j}^{\prime}
\big\} \cr
& & \hspace{50mm} (1 \le j \le m-1) \cr
\widetilde{f}_m^{\prime} &:=& \frac{1}{2}
\big\{
-Z_{2m-2, 2m-1}^{\prime}
-iZ_{2m-2}^{\prime}
\big\} 
\end{array}
\right. \cr
& & \left\{
\begin{array}{cclcl}
\widetilde{h}_j &:=& \frac{i}{2}
\big\{
Z_{2j-1}-Z_{2j+1}
\big\} 
& \qquad &(1 \le j \le m-1) \cr
\widetilde{h}_m &:=& 
iZ_{2m-1}
& & \cr
\widetilde{h}_j^{\prime} &:=& \frac{i}{2}
\big\{
Z_{2j-1}^{\prime}
-Z_{2j+1}^{\prime}
\big\} 
& \qquad &(1 \le j \le m-1) \cr
\widetilde{h}_m^{\prime} &:=& 
iZ_{2m-1}^{\prime}
& & 
\end{array}
\right. \cr
& & \left\{
\begin{array}{ccl}
\widetilde{e}_0 &:=& \frac{1}{2}
\big\{
i(X_{\alpha_1}-X_{-\alpha_1})+H_{\alpha_1}
\big\} \otimes t \cr
\widetilde{f}_0 &:=& \frac{1}{2}
\big\{
-i(X_{\alpha_1}-X_{-\alpha_1})+H_{\alpha_1}
\big\} \otimes t^{-1} \cr
\widetilde{h}_0 &:=& 
-iY_{\alpha_1}+\frac{K}{2}.
\end{array}
\right.
\label{eqn:(4.1.2)}
\end{eqnarray}

Then, by Lemma \ref{lemma:(4.1.1)}, 
one can easily check that these elements satisfy
the conditions of Chevalley generators for the following 
Dynkin diagrams according as $n=2m$ or $n=2m+1$,
letting $\widetilde{e}_j$ (resp. $\widetilde{e}_j^{\prime}$) be 
a root vector of a simple root $\widetilde{\alpha}_j$ (resp. 
$\widetilde{\alpha}_j^{\prime}$):

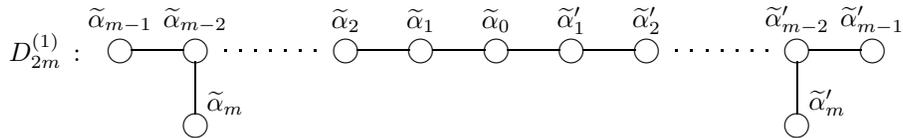
\begin{center}
\setlength{\unitlength}{1mm}
\begin{picture}(118,-18)
\put(15,-4){\circle{3}}
\put(25,-4){\circle{3}}
\put(45,-4){\circle{3}}
\put(55,-4){\circle{3}}
\put(65,-4){\circle{3}}
\put(75,-4){\circle{3}}
\put(85,-4){\circle{3}}
\put(105,-4){\circle{3}}
\put(115,-4){\circle{3}}
\put(25,-14){\circle{3}}
\put(105,-14){\circle{3}}
\put(16.5,-4){\line(1,0){7}}
\put(46.5,-4){\line(1,0){7}}
\put(56.5,-4){\line(1,0){7}}
\put(66.5,-4){\line(1,0){7}}
\put(76.5,-4){\line(1,0){7}}
\put(106.5,-4){\line(1,0){7}}
\put(25,-5.5){\line(0,-1){7}}
\put(105,-5.5){\line(0,-1){7}}
\put(15,0){\makebox(0,0){$\widetilde{\alpha}_{m-1}$}}
\put(25,0){\makebox(0,0){$\widetilde{\alpha}_{m-2}$}}
\put(45,0){\makebox(0,0){$\widetilde{\alpha}_{2}$}}
\put(55,0){\makebox(0,0){$\widetilde{\alpha}_{1}$}}
\put(65,0){\makebox(0,0){$\widetilde{\alpha}_{0}$}}
\put(75,0){\makebox(0,0){$\widetilde{\alpha}_{1}^{\prime}$}}
\put(85,0){\makebox(0,0){$\widetilde{\alpha}_{2}^{\prime}$}}
\put(105,0){\makebox(0,0){$\widetilde{\alpha}_{m-2}^{\prime}$}}
\put(115,0){\makebox(0,0){$\widetilde{\alpha}_{m-1}^{\prime}$}}
\put(109,-11){\makebox(0,0){$\widetilde{\alpha}_{m}^{\prime}$}}
\put(29,-11){\makebox(0,0){$\widetilde{\alpha}_{m}$}}
\put(5,-4){\makebox(0,0){$D_{2m}^{(1)} \, : $}}
\thicklines
\dottedline{2}(29,-4)(41,-4)
\dottedline{2}(89,-4)(101,-4)
\end{picture}
\end{center}

\vspace{15mm}

\begin{center}
\setlength{\unitlength}{1mm}
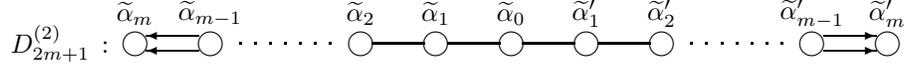
\begin{picture}(118,-8)
\put(15,-4){\circle{3}}
\put(25,-4){\circle{3}}
\put(45,-4){\circle{3}}
\put(55,-4){\circle{3}}
\put(65,-4){\circle{3}}
\put(75,-4){\circle{3}}
\put(85,-4){\circle{3}}
\put(105,-4){\circle{3}}
\put(115,-4){\circle{3}}
\put(23.5,-3){\vector(-1,0){7}}
\put(23.5,-5){\vector(-1,0){7}}
\put(46.5,-4){\line(1,0){7}}
\put(56.5,-4){\line(1,0){7}}
\put(66.5,-4){\line(1,0){7}}
\put(76.5,-4){\line(1,0){7}}
\put(106.5,-3){\vector(1,0){7}}
\put(106.5,-5){\vector(1,0){7}}
\put(15,0){\makebox(0,0){$\widetilde{\alpha}_{m}$}}
\put(25,0){\makebox(0,0){$\widetilde{\alpha}_{m-1}$}}
\put(45,0){\makebox(0,0){$\widetilde{\alpha}_{2}$}}
\put(55,0){\makebox(0,0){$\widetilde{\alpha}_{1}$}}
\put(65,0){\makebox(0,0){$\widetilde{\alpha}_{0}$}}
\put(75,0){\makebox(0,0){$\widetilde{\alpha}_{1}^{\prime}$}}
\put(85,0){\makebox(0,0){$\widetilde{\alpha}_{2}^{\prime}$}}
\put(105,0){\makebox(0,0){$\widetilde{\alpha}_{m-1}^{\prime}$}}
\put(115,0){\makebox(0,0){$\widetilde{\alpha}_{m}^{\prime}$}}
\put(5,-4){\makebox(0,0){$D_{2m+1}^{(2)} \, : $}}
\thicklines
\dottedline{2}(29,-4)(41,-4)
\dottedline{2}(89,-4)(101,-4)
\end{picture}
\end{center}

\vspace{8mm}

We denote by $\widetilde{\Lambda}_j$, 
$\widetilde{\Lambda}_j^{\prime}$ the fundamental weights 
corresponding to the simple coroot system 
$\widetilde{h}_j$, $\widetilde{h}_j^{\prime}$.
The action of $\widetilde{h}_j$ and $\widetilde{h}_j^{\prime}$ on 
the basis $v(c_1, \cdots, c_n)$ of $\ccc\{Q/2Q\}$ is calculated 
easily from Proposition \ref{prop:(4.1.1)} as follows:

\begin{prop}
\label{prop:(4.1.2)}
Let $v:=v(c_1, \cdots , c_n)$ where $c_1, \cdots, c_n \in
\{\pm 1\}$.  Then
\begin{enumerate}
\item[{\rm 1)}] \quad In the case $n=2m$,
$$
\begin{array}{ccll}
4\widetilde{h}_{j}v &=& c_{2j-1}
(1-c_{2j-1}c_{2j+1})
(1-c_{2m-1}c_{2m})v
& (1 \le j \le m-1) \cr
4\widetilde{h}_{m}v &=& c_{2m-3}
(1+c_{2m-3}c_{2m-1})
(1-c_{2m-1}c_{2m})v
& \cr
4\widetilde{h}_{j}^{\prime}v &=& c_{2j-1}
(1-c_{2j-1}c_{2j+1})
(1+c_{2m-1}c_{2m})v
& (1 \le j \le m-1) \cr
4\widetilde{h}_{m}^{\prime}v &=& c_{2m-3}
(1+c_{2m-3}c_{2m-1})
(1+c_{2m-1}c_{2m})v
& \cr
2\widetilde{h}_{0}v &=& (1-c_1)v \, . &
\end{array}
$$
\item[{\rm 2)}] \quad In the case $n=2m+1$,
$$
\begin{array}{ccll}
4\widetilde{h}_{j}v &=& c_{2j-1}
(1-c_{2j-1}c_{2j+1})(1-c_{2m}c_{2m+1})v
& (1 \le j \le m-1) \cr
2\widetilde{h}_{m}v &=& c_{2m-1}(1-c_{2m}c_{2m+1})v
& \cr
4\widetilde{h}_{j}^{\prime}v &=& c_{2j-1}
(1-c_{2j-1}c_{2j+1})(1+c_{2m}c_{2m+1})v
& (1 \le j \le m-1) \cr
2\widetilde{h}_{m}^{\prime}v &=& c_{2m-1}(1+c_{2m}c_{2m+1})v
& \cr
2\widetilde{h}_{0}v &=& (1-c_1)v \, . &
\end{array}
$$
\end{enumerate}
\end{prop}

Notice that, for a level one representation of a simply-laced 
algebra, a singular vector is characterized as an eigenvector of 
all $\widetilde{h}_j$'s and $\widetilde{h}_j^{\prime}$'s with 
non-negative integral eigenvalues.  From this, singular vectors are 
easily obtained by Proposition \ref{prop:(4.1.2)}.
And then, all other elements belonging to the invariant subspace 
spanned by a singular vector are obtained from 
Proposition \ref{prop:(4.1.1)}. The calculation is straightforward
and the result is stated as follows:

\begin{thm}
\label{thm:(4.1.1)}
Let $v=v(c_1, \cdots, c_n)$ where $c_1, \cdots , c_n \in \{\pm 1\}$. 
\begin{enumerate}
\item[{\rm 1)}] \quad In the case $n=2m$, 
\begin{enumerate}
\item[{\rm (i)}] \,\ $v$ is a singular 
vector if and only if $c_{2j-1}=1$ for $j=1, 2, \cdots, m-1$. 
And then the weight of $v$ is determined by $(c_{2m-1}, c_{2m})$ as 
is shown in the following table 
$$
\begin{array}{c|c}
\text{singular vector} & weight \cr
\noalign{\hrule height0.8pt}
\hfil
v \left(
\begin{array}{rcl}
1, c_2, 1, c_4 , \cdots , 1, & c_{2m-2}, & 1 \cr
& 1 &
\end{array}
\right)
& \widetilde{\Lambda}_m^{\prime} \\[2mm]
\hline
v \left(
\begin{array}{rcl}
1, c_2, 1, c_4 , \cdots , 1, & c_{2m-2}, & 1 \cr
& -1 &
\end{array}
\right)
& \widetilde{\Lambda}_m \\[2mm]
\hline
v \left(
\begin{array}{rcl}
1, c_2, 1, c_4 , \cdots , 1, & c_{2m-2}, & -1 \cr
& 1 &
\end{array}
\right)
& \widetilde{\Lambda}_{m-1}^{\prime} \\[2mm]
\hline
v \left(
\begin{array}{rcl}
1, c_2, 1, c_4 , \cdots , 1, & c_{2m-2}, & -1 \cr
& -1 &
\end{array}
\right)
& \widetilde{\Lambda}_{m-1} \\[2mm]
\noalign{\hrule height0.8pt}
\end{array}
$$
for any choice of $c_2, c_4, \cdots, c_{2m-2} \in \{\pm 1\}$.

\item[{\rm (ii)}] \,\ Given a singular vector 
$$v_0 \, := \, v
\left(
\begin{array}{rcl}
1, c_2, 1, c_4 , \cdots , 1, & c_{2m-2}, & c_{2m-1} \cr
& c_{2m} &
\end{array}
\right),$$
the $\ccc$-linear span of all elements
$$v
\left(
\begin{array}{rcl}
b_1, c_2, b_3, c_4 , \cdots , b_{2m-3}, & c_{2m-2}, & b_{2m-1} \cr
& b_{2m} &
\end{array}
\right)$$
satisfying the conditions
\begin{enumerate}
\item[{\rm (a)}] \qquad ${\displaystyle
b_j \in \{\pm 1\} }$ \quad for all $j$,
\item[{\rm (b)}] \qquad $b_{2m-1}b_{2m} =c_{2m-1}c_{2m}$,
\item[{\rm (c)}] \qquad ${\displaystyle
\prod_{\substack{1 \le j \le 2m-1 \cr
j= \, {\rm odd}}} b_j = c_{2m-1},}$
\end{enumerate}
tensored with $\ccc [ x^{(j)}_r ; 1 \le j \le 2m , \,\ 
r \in \nnn_{\rm odd} ]$ is the irreducible $D^{(1)}_{2m}$-module
with the highest weight vector $v_0$.
\end{enumerate}

\item[{\rm 2)}] \quad In the case $n=2m+1$, 
\begin{enumerate}
\item[{\rm (i)}] \,\ $v$ is a singular 
vector if and only if $c_{2j-1}=1$ for $j=1, 2, \cdots, m$. 
And then the weight of $v$ is determined by $(c_{2m}, c_{2m+1})$ as 
is shown in the following table:
$$
\begin{array}{c|c}
\text{singular vector} & weight \cr
\noalign{\hrule height0.8pt}
\hfil
v \left(
\begin{array}{rcl}
1, c_2, 1, c_4 , \cdots , 1 , \, c_{2m-2}, & 1, & \pm 1 \cr
& \pm 1 &
\end{array}
\right)
& \widetilde{\Lambda}_m^{\prime} \\[2mm]
\hline
v \left(
\begin{array}{rcl}
1, c_2, 1, c_4 , \cdots , 1, \, c_{2m-2}, & 1, & \pm 1 \cr
& \mp 1 &
\end{array}
\right)
& \widetilde{\Lambda}_m \\[2mm]
\noalign{\hrule height0.8pt}
\end{array}
$$
for any choice of $c_2, c_4, \cdots, c_{2m-2} \in \{\pm 1\}$.

\item[{\rm (ii)}] \,\ Given a singular vector 
$$v_0 \, := \, v
\left(
\begin{array}{rcl}
1, c_2, 1, c_4 ,1,  \cdots , 1, \, c_{2m-2}, & 1, & c_{2m} \cr
& c_{2m+1} &
\end{array}
\right),$$
the $\ccc$-linear span of all elements
$$v
\left(
\begin{array}{rcl}
b_1, c_2, b_3, c_4 , \cdots , b_{2m-3}, c_{2m-2}, 
& b_{2m-1}, & c_{2m} \cr
& c_{2m+1} &
\end{array}
\right),$$
where $b_j \in \{\pm 1\}$ for all $j$,
tensored with $\ccc [ x^{(j)}_r ; 1 \le j \le 2m+1 , \,\ 
r \in \nnn_{\rm odd} ]$ is the irreducible $D^{(2)}_{2m+1}$-module
with the highest weight vector $v_0$.
\end{enumerate}
\end{enumerate}
\end{thm}

In the case $D_n$, one can choose an orthonormal basis 
$\varepsilon_1, \cdots, \varepsilon_n$ of $\hhh^{\ast}$ such that
$$\alpha_j = \varepsilon_j-\varepsilon_{j+1} \,\ (1 \le j \le n-1) 
\qquad \text{and} \qquad 
\alpha_n =\varepsilon_{n-1}+\varepsilon_n,$$ 
and let $S_1, \cdots, S_n$ be its dual basis of $\hhh$.
Then the set of positive roots is given by 
$$\Delta_+ \,\ = \,\ \{ \varepsilon_j \pm \varepsilon_k \,\ ; \,\ 
1 \le j < k \le n \},$$
and, for each $\varepsilon_j \pm \varepsilon_k$, the 
operator $\Gamma_{\varepsilon_j \pm \varepsilon_k}(z)$ is 
written as follows:
$$
\Gamma_{\varepsilon_j \pm \varepsilon_k}(z)
= 
\frac{1}{2}e^{\varepsilon_j \pm \varepsilon_k}
\nu(\varepsilon_j \pm \varepsilon_k, \,\ \cdot \,\ )
\widetilde{U}_{\varepsilon_j \pm \varepsilon_k}(z)
$$
where
\begin{eqnarray*}
\widetilde{U}_{\varepsilon_j \pm \varepsilon_k}(z)
&:=&
\exp \left( \sqrt{2}
\sum_{r \in \nnn_{\rm odd}}
\big(x^{(j)}_r \pm x^{(k)}_r\big)z^r\right)
\cr
&\times&
\exp \left( -\sqrt{2} \sum_{r \in \nnn_{\rm odd}}
\left(
\frac{\partial}{\partial x^{(j)}_r} 
\pm 
\frac{\partial}{\partial x^{(k)}_r} 
\right)
\frac{z^{-r}}{r}
\right).
\end{eqnarray*}

\begin{ex}
\label{ex:(4.1.1)}
We consider $D_4$ with the oriented Dynkin diagram
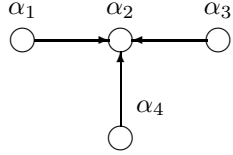
\begin{center}
\setlength{\unitlength}{1mm}
\begin{picture}(30,-18)
\put(1,-4){\circle{3}}
\put(14,-4){\circle{3}}
\put(27,-4){\circle{3}}
\put(14,-17){\circle{3}}
\put(2.5,-4){\vector(1,0){10}}
\put(25.5,-4){\vector(-1,0){10}}
\put(14,-15.5){\vector(0,1){10}}
\put(1,0){\makebox(0,0){$\alpha_1$}}
\put(14,0){\makebox(0,0){$\alpha_2$}}
\put(27,0){\makebox(0,0){$\alpha_3$}}
\put(18,-13){\makebox(0,0){$\alpha_4$}}
\end{picture}
\end{center}

\vspace{18mm}

\noindent
We write simply \, 
$\big(
\begin{array}{c}
m_1 \, m_2 \,  m_3 \\[-1mm]
m_4
\end{array}
\big)$ \, 
for $\alpha=\sum^4_{j=1}m_j\alpha_j$. 
Fix $c_1, c_2, c_3, c_4 \in \{\pm 1 \}$, and put
$$
v:=v(c_1, c_2, c_3, c_4) \qquad \text{and} \qquad 
v^{\prime}:=v(-c_1, c_2, -c_3, -c_4). $$

Then, by Lemma \ref{lemma:(3.1)}, one sees that the $\ccc$-linear 
span of $v$ and $v^{\prime}$ is invariant under the action 
of $\widehat{X}_{\alpha_j}$'s and so $\widehat{X}_{\alpha}$
for all $\alpha \in \Delta_+$. 
To write down the action of $\widehat{X}_{\alpha}$ on 
the space $\ccc v \oplus \ccc v^{\prime}$ explicitly, we recall 
the Pauli's spin matrices 
$$
\sigma_1 :=
\begin{pmatrix}
0 & 1 \cr
1 & 0
\end{pmatrix}
\qquad 
\sigma_2 :=
\begin{pmatrix}
0 & -i \cr
i & 0
\end{pmatrix}
\qquad 
\sigma_3 :=
\begin{pmatrix}
1 & 0 \cr
0 & -1
\end{pmatrix}
$$
and divide the set of positive roots into the disjoint union
of three parts:
$$
\begin{array}{ccc}
\Delta^{(1)}_+ &:=& \{
\varepsilon_1 \pm \varepsilon_4 , \,\ 
\varepsilon_2 \pm \varepsilon_3 \},
\cr
\Delta^{(2)}_+ &:=& \{
\varepsilon_1 \pm \varepsilon_3 , \,\ 
\varepsilon_2 \pm \varepsilon_4 \},
\cr
\Delta^{(3)}_+ &:=& \{
\varepsilon_1 \pm \varepsilon_2 , \,\ 
\varepsilon_3 \pm \varepsilon_4 \}.
\end{array}
$$
Then this decomposition \, $\Delta_+ =
\Delta^{(1)}_+ \cup \Delta^{(2)}_+ \cup \Delta^{(3)}_+$ \, 
has the following properties:
\begin{enumerate}
\item[{\rm (i)}] \quad $
\alpha, \beta \in \Delta^{(i)} \,\ \Longrightarrow \,\ 
\alpha \pm \beta \not\in \Delta $.
\item[{\rm (ii)}] \quad For $\{i, j, k \}=\{1, 2,3\}$, 
$$
\alpha \in \Delta^{(i)}_+, \,\ \beta \in \Delta^{(j)}_+ \,\
\text{and} \,\ \alpha+\beta \in \Delta
\quad \Longrightarrow \quad 
\alpha+\beta \in \Delta^{(k)}_+.$$
\end{enumerate}

For $\alpha \in \Delta_+$, 
$\widehat{X}_{\alpha}v$ and $\widehat{X}_{\alpha}v^{\prime}$ are 
computed by using Lemma \ref{lemma:(3.1)} and \eqref{eqn:(3.4)}, 
and are explicitly written as follows:
\begin{eqnarray*}
\Delta^{(1)}_+
&: & \left\{
\begin{array}{ccccl}
2\widehat{X}_{\left(
\begin{subarray}{c}
010 \cr
0
\end{subarray}
\right) }
&=&
2\widehat{X}_{\varepsilon_2-\varepsilon_3}
&=&
ic_2\sigma_1 \cr
2\widehat{X}_{\left(
\begin{subarray}{c}
111 \cr
0
\end{subarray} \right) }
&=&
2\widehat{X}_{\varepsilon_1-\varepsilon_4}
&=&
ic_1c_2c_3\sigma_1 \cr
2\widehat{X}_{\left(
\begin{subarray}{c}
110 \cr
1
\end{subarray} \right) }
&=&
2\widehat{X}_{\varepsilon_1+\varepsilon_4}
&=&
ic_1c_2c_4\sigma_1 \cr
2\widehat{X}_{\left(
\begin{subarray}{c}
011 \cr
1
\end{subarray} \right) }
&=&
2\widehat{X}_{\varepsilon_2+\varepsilon_3}
&=&
ic_2c_3c_4\sigma_1 
\end{array}
\right. 
\cr
\Delta^{(2)}_+
&: & \left\{
\begin{array}{ccccl}
2\widehat{X}_{\left(
\begin{subarray}{c}
110 \cr
0
\end{subarray} \right) }
&=&
2\widehat{X}_{\varepsilon_1-\varepsilon_3}
&=&
-ic_1c_2\sigma_2 \cr
2\widehat{X}_{\left(
\begin{subarray}{c}
011 \cr
0
\end{subarray} \right) }
&=&
2\widehat{X}_{\varepsilon_2-\varepsilon_4}
&=&
-ic_2c_3\sigma_2 \cr
2\widehat{X}_{\left(
\begin{subarray}{c}
010 \cr
1
\end{subarray} \right) }
&=&
2\widehat{X}_{\varepsilon_2+\varepsilon_4}
&=&
-ic_2c_4\sigma_2 \cr
2\widehat{X}_{\left(
\begin{subarray}{c}
111 \cr
1
\end{subarray} \right) }
&=&
2\widehat{X}_{\varepsilon_1+\varepsilon_3}
&=&
-ic_1c_2c_3c_4\sigma_2 \cr
\end{array}
\right. \cr
\Delta^{(3)}_+
&: & \left\{
\begin{array}{ccccl}
2\widehat{X}_{\left(
\begin{subarray}{c}
100 \cr
0
\end{subarray} \right) }
&=&
2\widehat{X}_{\varepsilon_1-\varepsilon_2}
&=&
-ic_1\sigma_3 \cr
2\widehat{X}_{\left(
\begin{subarray}{c}
001 \cr
0
\end{subarray} \right) }
&=&
2\widehat{X}_{\varepsilon_3-\varepsilon_4}
&=&
-ic_3\sigma_3 \cr
2\widehat{X}_{\left(
\begin{subarray}{c}
000 \cr
1
\end{subarray} \right) }
&=&
2\widehat{X}_{\varepsilon_3+\varepsilon_4}
&=&
-ic_4\sigma_3 \cr
2\widehat{X}_{\left(
\begin{subarray}{c}
121 \cr
1
\end{subarray} \right) }
&=&
2\widehat{X}_{\varepsilon_1+\varepsilon_2}
&=&
-ic_1c_3c_4\sigma_3 \, .
\end{array}
\right.
\end{eqnarray*}
Then, in paticular letting $(c_1, c_2, c_3, c_4):=
(-1, 1 ,-1, -1)$, one obtains the irreducible 
representation $\pi$ of $D^{(1)}_4$ on the space 
\begin{eqnarray*}
\widetilde{V} & \hspace{-3mm} :=& \hspace{-3mm}
\left(v \otimes 
\ccc [ x^{(j)}_r ; 1 \le j \le 4, \, r \in \nnn_{\rm odd} ]
\right)
\oplus \left(
v^{\prime} \otimes 
\ccc [ x^{(j)}_r ; 1 \le j \le 4, \, r \in \nnn_{\rm odd} ]
 \right)
\cr
& \hspace{-3mm} \cong &  \hspace{-3mm}
\ccc [ x^{(j)}_r ; 1 \le j \le 4, \, r \in \nnn_{\rm odd} ]
\oplus 
\ccc [ x^{(j)}_r ; 1 \le j \le 4, \, r \in \nnn_{\rm odd} ]
\end{eqnarray*}
as follows:
$$
\pi \,\ : \,\ \left\{
\begin{array}{ccl}
S_j \otimes t^{r} &\longmapsto&
{\displaystyle
\frac{1}{\sqrt{2}} \cdot 
\frac{\partial}{\partial x^{(j)}_r} \cdot \sigma_0} \cr
S_j \otimes t^{-r} &\longmapsto&
{\displaystyle
\frac{r}{\sqrt{2}} \cdot x^{(j)}_r \sigma_0 }
\end{array}
\qquad (j=1,2,3,4 ; \,\ r \in \nnn_{\rm odd})
\right.
$$
and
$$
\pi \,\ : \,\ \left\{
\begin{array}{ccl}
\widetilde{X}_{\varepsilon_j \pm \varepsilon_k}(z) 
&\longmapsto &{\displaystyle 
\frac{i}{2} 
\widetilde{U}_{\varepsilon_j \pm \varepsilon_k}(z) 
\cdot \sigma_p} \qquad \quad (\text{if} \,\ 
\varepsilon_j \pm \varepsilon_k \in \Delta^{(p)}_+
),
\cr
d &\longmapsto & {\displaystyle - \sum^4_{j=1}
\sum_{r \in \nnn_{\rm odd}}
rx^{(j)}_r \frac{\partial}{\partial x^{(j)}_r}} \cr
K & \longmapsto & 1 \,\ (= \,\ \text{the identity operator})
\end{array}
\right.
$$
where we put $\sigma_0 := 
\begin{pmatrix}
1 & 0 \cr
0 & 1
\end{pmatrix}
.$
\end{ex}

\subsection{The case $A_n$}

For $A_n$, we consider the following orientation of Dynkin
diagram according as $n$ is odd or even:
\begin{center}
\setlength{\unitlength}{1mm}
\begin{picture}(103,-8)
\put(15,-4){\circle{3}}
\put(25,-4){\circle{3}}
\put(35,-4){\circle{3}}
\put(45,-4){\circle{3}}
\put(55,-4){\circle{3}}
\put(75,-4){\circle{3}}
\put(87,-4){\circle{3}}
\put(99,-4){\circle{3}}
\put(16.5,-4){\vector(1,0){7}}
\put(36.5,-4){\vector(1,0){7}}
\put(76.5,-4){\vector(1,0){9}}
\put(33.5,-4){\vector(-1,0){7}}
\put(53.5,-4){\vector(-1,0){7}}
\put(97.5,-4){\vector(-1,0){9}}
\put(15,0){\makebox(0,0){$\alpha_1$}}
\put(25,0){\makebox(0,0){$\alpha_2$}}
\put(35,0){\makebox(0,0){$\alpha_3$}}
\put(45,0){\makebox(0,0){$\alpha_4$}}
\put(55,0){\makebox(0,0){$\alpha_5$}}
\put(75,0){\makebox(0,0){$\alpha_{2m-3}$}}
\put(87,0){\makebox(0,0){$\alpha_{2m-2}$}}
\put(99,0){\makebox(0,0){$\alpha_{2m-1}$}}
\put(0,-4){\makebox(0,0){$A_{2m-1} \,\ : $}}
\thicklines
\dottedline{2}(59,-4)(71,-4)
\end{picture}
\end{center}

\vspace{5mm}

\begin{center}
\setlength{\unitlength}{1mm}
\begin{picture}(103,-8)
\put(15,-4){\circle{3}}
\put(25,-4){\circle{3}}
\put(35,-4){\circle{3}}
\put(45,-4){\circle{3}}
\put(55,-4){\circle{3}}
\put(75,-4){\circle{3}}
\put(87,-4){\circle{3}}
\put(99,-4){\circle{3}}
\put(16.5,-4){\vector(1,0){7}}
\put(36.5,-4){\vector(1,0){7}}
\put(88.5,-4){\vector(1,0){9}}
\put(33.5,-4){\vector(-1,0){7}}
\put(53.5,-4){\vector(-1,0){7}}
\put(85.5,-4){\vector(-1,0){9}}
\put(15,0){\makebox(0,0){$\alpha_1$}}
\put(25,0){\makebox(0,0){$\alpha_2$}}
\put(35,0){\makebox(0,0){$\alpha_3$}}
\put(45,0){\makebox(0,0){$\alpha_4$}}
\put(55,0){\makebox(0,0){$\alpha_5$}}
\put(75,0){\makebox(0,0){$\alpha_{2m-2}$}}
\put(87,0){\makebox(0,0){$\alpha_{2m-1}$}}
\put(99,0){\makebox(0,0){$\alpha_{2m}$}}
\put(0,-4){\makebox(0,0){$A_{2m} \,\ : $}}
\thicklines
\dottedline{2}(59,-4)(71,-4)
\end{picture}
\end{center}

\vspace{8mm}

\begin{lemma} 
\label{lemma:(4.2.1)}
Let $n=2m$ or $2m-1$ in $A_n$. Then the action 
of $Y_{\alpha_j}$ to an element $v:=v(c_1, \cdots , c_n)$ 
 $(c_1, \cdots, c_n \in \{\pm 1\})$ is given as follows:
\begin{enumerate}
\item[{\rm 1)}] \quad $
2Y_{\alpha_{2j-1}}v(c_1, \cdots, c_n)
=
-ic_{2j-1}v(c_1, \cdots, c_n)  \qquad (1 \le j \le m),
$
\item[{\rm 2)}] \quad $
2Y_{\alpha_{2j}}v(c_1, \cdots , c_n)$
$$= \left\{
\begin{array}{lcl}
-ic_{2j}v(c_1, \cdots , -c_{2j-1}, c_{2j}, -c_{2j+1}, \, \cdots, c_n) 
& &\text{if} \,\ 2j<n \cr
-ic_{2m}v(c_1, \cdots , c_{2m-2}, -c_{2m-1}, c_{2m})
& &\text{if} \,\ 2j=n=2m. 
\end{array}
\right.
$$
\end{enumerate}
\end{lemma}

For $\ggg=A_n=sl(n+1, \ccc)$, the $\sigma$-fixed subalgebra 
$\ggg_{\bar{0}}$ is $so(n+1, \ccc)$, which is a simple Lie algebra 
of type $D_m$ if $n=2m-1$ and $B_m$ if $n=2m$.
For $1 \le j \le k \le n$ we put
$$Y_{j,k} := Y_{\alpha_j + \cdots +\alpha_k} \qquad 
\text{and} \qquad Y_j :=Y_{j,j} =Y_{\alpha_j},$$
and define elements $\widetilde{e}_j, \widetilde{f}_j, 
\widetilde{h}_j$ ($0 \le j \le m$) in $\widehat{\ggg}(\sigma)$ as 
follows:

\noindent
In the case $n=2m-1$:
\begin{eqnarray}
& & \left\{
\begin{array}{ccl}
\widetilde{e}_j &:=& \frac{1}{2}
\big\{
Y_{2j-1, 2j}-Y_{2j, 2j+1}-iY_{2j-1, 2j+1}-iY_{2j}
\big\} \cr
& & \hspace{50mm} (1 \le j \le m-1) \cr
\widetilde{e}_m &:=& \frac{1}{2}
\big\{
Y_{2m-3, 2m-2}+Y_{2m-2, 2m-1}+iY_{2m-3, 2m-1}-iY_{2m-2}
\big\}
\end{array}
\right. \cr
& & \left\{
\begin{array}{ccl}
\widetilde{f}_j &:=& \frac{1}{2}
\big\{
-Y_{2j-1, 2j}+Y_{2j, 2j+1}-iY_{2j-1, 2j+1}-iY_{2j} 
\big\} \cr
& & \hspace{50mm} (1 \le j \le m-1) \cr
\widetilde{f}_m &:=& \frac{1}{2}
\big\{
-Y_{2m-3, 2m-2}-Y_{2m-2, 2m-1}+iY_{2m-3, 2m-1}-iY_{2m-2}
\big\} 
\end{array}
\right. \cr
& & \left\{
\begin{array}{cclcl}
\widetilde{h}_j &:=& i\big\{ Y_{2j-1}-Y_{2j+1} \big\} 
& \qquad \qquad &(1 \le j \le m-1) \cr
\widetilde{h}_m &:=& i\big\{ Y_{2m-3}+Y_{2m-1} \big\} 
& & 
\end{array}
\right. \cr
& & \left\{
\begin{array}{ccl}
\widetilde{e}_0 &:=& \frac{1}{2}
\big\{ i(X_{\alpha_1}-X_{-\alpha_1})+H_{\alpha_1} \big\} 
\otimes t \cr
\widetilde{f}_0 &:=& \frac{1}{2}
\big\{ -i(X_{\alpha_1}-X_{-\alpha_1})+H_{\alpha_1} \big\} 
\otimes t^{-1} \cr
\widetilde{h}_0 &:=& -iY_{\alpha_1}+\frac{K}{2}.
\end{array}
\right.
\label{eqn:(4.2.1)}
\end{eqnarray}
In the case $n=2m$:
\begin{eqnarray}
& & \left\{
\begin{array}{ccl}
\widetilde{e}_j &:=& \frac{1}{2}
\big\{
Y_{2j-1, 2j}-Y_{2j, 2j+1}-iY_{2j-1, 2j+1}-iY_{2j}
\big\} \cr
& & \hspace{50mm} (1 \le j \le m-1) \cr
\widetilde{e}_m &:=& Y_{2m-2, 2m-1}-iY_{2m-2} 
\end{array} \right. \cr
& & \left\{
\begin{array}{ccl}
\widetilde{f}_j &:=& \frac{1}{2}
\big\{
-Y_{2j-1, 2j}+Y_{2j, 2j+1}-iY_{2j-1, 2j+1}-iY_{2j}
\big\} \cr
& & \hspace{50mm} (1 \le j \le m-1) \cr
\widetilde{f}_m &:=& -Y_{2m-2, 2m-1}-iY_{2m-2}
\end{array}
\right. \cr
& & \left\{
\begin{array}{cclcl}
\widetilde{h}_j &:=& i
\big\{ Y_{2j-1}-Y_{2j+1} \big\} 
& \qquad \qquad &(1 \le j \le m-1) \cr
\widetilde{h}_m &:=& 2iY_{2m-1} & & 
\end{array} 
\right. \cr
& & \left\{
\begin{array}{ccl}
\widetilde{e}_0 &:=& \frac{1}{2} \big\{
i(X_{\alpha_1}-X_{-\alpha_1})+H_{\alpha_1}
\big\} \otimes t \cr
\widetilde{f}_0 &:=& \frac{1}{2} \big\{
-i(X_{\alpha_1}-X_{-\alpha_1})+H_{\alpha_1} 
\big\} \otimes t^{-1} \cr
\widetilde{h}_0 &:=& -iY_{\alpha_1}+\frac{K}{2}.
\end{array}
\right.
\label{eqn:(4.2.2)}
\end{eqnarray}

Then one can easily check that these elements satisfy the 
conditions of Chevalley generators for the following Dynkin
diagrams according as $n=2m-1$ or $n=2m$, letting 
$\widetilde{e}_j$ be a root vector of a simple root
$\widetilde{\alpha}_j$:

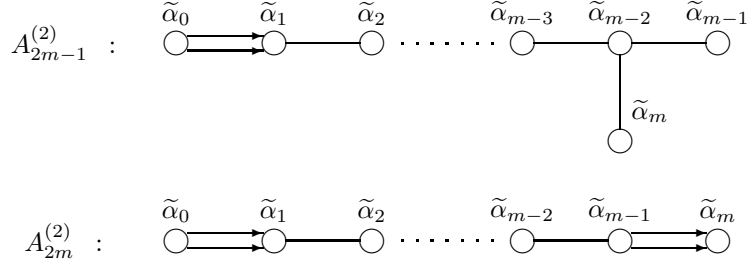
\begin{center}
\setlength{\unitlength}{1mm}
\begin{picture}(90,-19)
\put(15,-4){\circle{3}}
\put(28,-4){\circle{3}}
\put(41,-4){\circle{3}}
\put(61,-4){\circle{3}}
\put(74,-4){\circle{3}}
\put(87,-4){\circle{3}}
\put(74,-17){\circle{3}}
\put(16.5,-3){\vector(1,0){10}}
\put(16.5,-5){\vector(1,0){10}}
\put(29.5,-4){\line(1,0){10}}
\put(62.5,-4){\line(1,0){10}}
\put(75.5,-4){\line(1,0){10}}
\put(74,-5,5){\line(0,-1){10}}
\put(15,0){\makebox(0,0){$\widetilde{\alpha}_0$}}
\put(28,0){\makebox(0,0){$\widetilde{\alpha}_1$}}
\put(41,0){\makebox(0,0){$\widetilde{\alpha}_2$}}
\put(61,0){\makebox(0,0){$\widetilde{\alpha}_{m-3}$}}
\put(74,0){\makebox(0,0){$\widetilde{\alpha}_{m-2}$}}
\put(87,0){\makebox(0,0){$\widetilde{\alpha}_{m-1}$}}
\put(78,-13){\makebox(0,0){$\widetilde{\alpha}_m$}}
\put(0,-4){\makebox(0,0){$A^{(2)}_{2m-1} \,\ : $}}
\thicklines
\dottedline{2}(45,-4)(57,-4)
\end{picture}
\end{center}

\vspace{15mm}

\begin{center}
\setlength{\unitlength}{1mm}
\begin{picture}(90,-6)
\put(15,-4){\circle{3}}
\put(28,-4){\circle{3}}
\put(41,-4){\circle{3}}
\put(61,-4){\circle{3}}
\put(74,-4){\circle{3}}
\put(87,-4){\circle{3}}
\put(16.5,-3){\vector(1,0){10}}
\put(16.5,-5){\vector(1,0){10}}
\put(29.5,-4){\line(1,0){10}}
\put(62.5,-4){\line(1,0){10}}
\put(75.5,-3){\vector(1,0){10}}
\put(75.5,-5){\vector(1,0){10}}
\put(15,0){\makebox(0,0){$\widetilde{\alpha}_0$}}
\put(28,0){\makebox(0,0){$\widetilde{\alpha}_1$}}
\put(41,0){\makebox(0,0){$\widetilde{\alpha}_2$}}
\put(61,0){\makebox(0,0){$\widetilde{\alpha}_{m-2}$}}
\put(74,0){\makebox(0,0){$\widetilde{\alpha}_{m-1}$}}
\put(87,0){\makebox(0,0){$\widetilde{\alpha}_{m}$}}
\put(0,-4){\makebox(0,0){$A^{(2)}_{2m} \,\ : $}}
\thicklines
\dottedline{2}(45,-4)(57,-4)
\end{picture}
\end{center}

\vspace{7mm}

\noindent
The action of $\widetilde{h}_j$ on the basis 
$v(c_1, \cdots, c_n)$ of $\ccc\{Q/2Q\}$ is calculated from 
Lemma \ref{lemma:(4.2.1)} as follows:

\begin{prop}
\label{prop:(4.2.1)}
Let $v:=v(c_1, \cdots, c_n)$ where $c_1, \cdots, c_n \in 
\{\pm 1\}$.
\begin{enumerate}
\item[{\rm 1)}] \quad In the case $n=2m-1$,
$$
2\widetilde{h}_jv= \left\{
\begin{array}{lcl}
(1-c_1)v & & (j=0), \cr
c_{2j-1}(1-c_{2j-1}c_{2j+1})v & & (1 \le j \le m-1), \cr
c_{2m-3}(1+c_{2m-3}c_{2m-1})v & & (j = m). 
\end{array}
\right.
$$
\item[{\rm 2)}] \quad In the case $n=2m$,
$$
2\widetilde{h}_jv= \left\{
\begin{array}{lcl}
(1-c_1)v & & (j=0), \cr
c_{2j-1}(1-c_{2j-1}c_{2j+1})v & & (1 \le j \le m-1), \cr
2c_{2m-1}v & & (j = m). 
\end{array}
\right.
$$
\end{enumerate}
\end{prop}

From this proposition one obtains the following theorem, where 
$\widetilde{\Lambda}_j$'s are the fundamental integral form 
corresponding to the simple coroots system $\widetilde{h}_j$.

\begin{thm}
\label{thm:(4.2.1)}
Let $v:=v(c_1, \cdots, c_n)$ where $c_1, \cdots, c_n \in 
\{\pm 1\}$.
\begin{enumerate}
\item[{\rm 1)}] \quad In the case $n=2m-1$,
\begin{enumerate}
\item[{\rm (i)}] \,\ 
$v$ is a singular vector if and only if $c_{2j-1}=1$ for 
$j=1, 2, \cdots, m-1$. And then the weight of $v$ is determined
by $c_{2m-1}$ as follows:
$$
\begin{array}{ccc}
\text{singular vector} & & \text{weight} \cr
v(1, c_2, 1, c_4, 1, c_6, \cdots, 1, c_{2m-2}, 1) & : &
\widetilde{\Lambda}_m \cr
v(1, c_2, 1, c_4, 1, c_6, \cdots, 1, c_{2m-2}, -1) & : &
\widetilde{\Lambda}_{m-1} 
\end{array}
$$
for any choice of $c_2, c_4, \cdots, c_{2m-2} \in \{\pm 1\}$.
\item[{\rm (ii)}] \,\ 
Given a singular vector 
$$v_0 \, := \, 
v(1, c_2, 1, c_4, 1, c_6, \cdots, 1, c_{2m-2}, c_{2m-1}),$$
the $\ccc$-linear span of all elements 
$$v(b_1, c_2, b_3, c_4, b_5, c_6, \cdots, b_{2m-3}, 
c_{2m-2}, b_{2m-1})$$
satisfying the conditions
\begin{enumerate}
\item[{\rm (a)}] \qquad 
$b_j \in \{\pm 1\}$ \qquad for all \, $j$
\item[{\rm (b)}] \qquad ${\displaystyle
\prod_{\substack{1 \le j \le 2m-1 \cr
j= \, {\rm odd}}} b_j = c_{2m-1} }$
\end{enumerate}
tensored with $\ccc[ x^{(j)}_r ; 1 \le j \le 2m-1, \,\ r \in 
\nnn_{\rm odd} ] $ is the irreducible $A^{(2)}_{2m-1}$-module
with the highest weight vector $v_0$.
\end{enumerate}

\item[{\rm 2)}] \quad In the case $n=2m$,
\begin{enumerate}
\item[{\rm (i)}] \,\ 
$v$ is a singular vector if and only if $c_{2j-1}=1$ for 
$j=1, 2, \cdots, m$. And then the weight of 
$$v(1, c_2, 1, c_4, 1, c_6, \cdots, 1, c_{2m-2}, 1, c_{2m})$$
is $\widetilde{\Lambda}_m$ for any choice of $c_2, c_4, 
\cdots, c_{2m} \in \{\pm 1\}$.
\item[{\rm (ii)}] \,\ 
Given a singular vector 
$$v_0 \, := \, 
v(1, c_2, 1, c_4, 1, c_6, \cdots, 1, c_{2m-2}, 1, c_{2m}),$$
the $\ccc$-linear span of all elements 
$$v(b_1, c_2, b_3, c_4, b_5, c_6, \cdots, b_{2m-3}, 
c_{2m-2}, b_{2m-1}, c_{2m}),$$
where $b_j \in \{\pm 1\}$ for all $j$,
tensored with $\ccc[ x^{(j)}_r ; 1 \le j \le 2m, \,\ r \in 
\nnn_{\rm odd} ] $ is the irreducible $A^{(2)}_{2m}$-module
with the highest weight vector $v_0$.
\end{enumerate}
\end{enumerate}
\end{thm}

\subsection{The case $E_n$ ($n=6,7,8$)}

For $E_n$, we consider the following orientation of Dynkin diagrams:
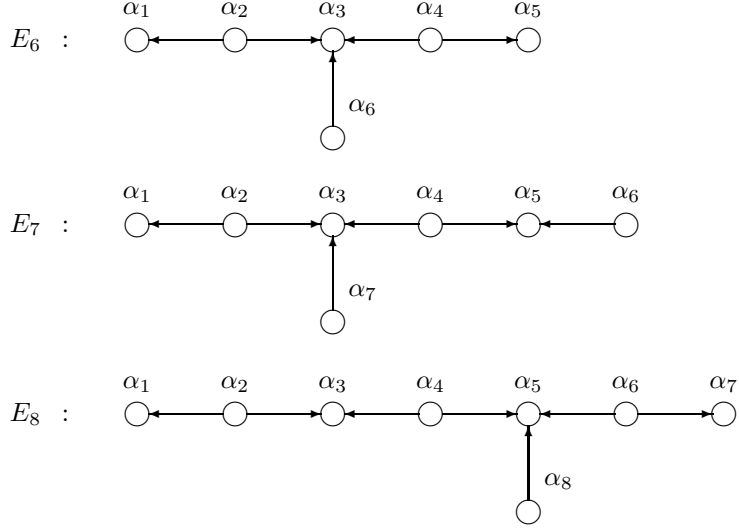
\begin{center}
\setlength{\unitlength}{1mm}
\begin{picture}(100,-19)
\put(18,-4){\circle{3}}
\put(31,-4){\circle{3}}
\put(44,-4){\circle{3}}
\put(57,-4){\circle{3}}
\put(70,-4){\circle{3}}
\put(44,-17){\circle{3}}
\put(32.5,-4){\vector(1,0){10}}
\put(58.5,-4){\vector(1,0){10}}
\put(29.5,-4){\vector(-1,0){10}}
\put(55.5,-4){\vector(-1,0){10}}
\put(44,-15.5){\vector(0,1){10}}
\put(18,0){\makebox(0,0){$\alpha_1$}}
\put(31,0){\makebox(0,0){$\alpha_2$}}
\put(44,0){\makebox(0,0){$\alpha_3$}}
\put(57,0){\makebox(0,0){$\alpha_4$}}
\put(70,0){\makebox(0,0){$\alpha_5$}}
\put(48,-13){\makebox(0,0){$\alpha_6$}}
\put(5,-4){\makebox(0,0){$E_6 \,\ :$}}
\end{picture}
\end{center}

\vspace{14mm}

\begin{center}
\setlength{\unitlength}{1mm}
\begin{picture}(100,-19)
\put(18,-4){\circle{3}}
\put(31,-4){\circle{3}}
\put(44,-4){\circle{3}}
\put(57,-4){\circle{3}}
\put(70,-4){\circle{3}}
\put(83,-4){\circle{3}}
\put(44,-17){\circle{3}}
\put(32.5,-4){\vector(1,0){10}}
\put(58.5,-4){\vector(1,0){10}}
\put(29.5,-4){\vector(-1,0){10}}
\put(55.5,-4){\vector(-1,0){10}}
\put(81.5,-4){\vector(-1,0){10}}
\put(44,-15.5){\vector(0,1){10}}
\put(18,0){\makebox(0,0){$\alpha_1$}}
\put(31,0){\makebox(0,0){$\alpha_2$}}
\put(44,0){\makebox(0,0){$\alpha_3$}}
\put(57,0){\makebox(0,0){$\alpha_4$}}
\put(70,0){\makebox(0,0){$\alpha_5$}}
\put(83,0){\makebox(0,0){$\alpha_6$}}
\put(48,-13){\makebox(0,0){$\alpha_7$}}
\put(5,-4){\makebox(0,0){$E_7 \,\ :$}}
\end{picture}
\end{center}

\vspace{14mm}

\begin{center}
\setlength{\unitlength}{1mm}
\begin{picture}(100,-19)
\put(18,-4){\circle{3}}
\put(31,-4){\circle{3}}
\put(44,-4){\circle{3}}
\put(57,-4){\circle{3}}
\put(70,-4){\circle{3}}
\put(83,-4){\circle{3}}
\put(96,-4){\circle{3}}
\put(70,-17){\circle{3}}
\put(32.5,-4){\vector(1,0){10}}
\put(58.5,-4){\vector(1,0){10}}
\put(84.5,-4){\vector(1,0){10}}
\put(29.5,-4){\vector(-1,0){10}}
\put(55.5,-4){\vector(-1,0){10}}
\put(81.5,-4){\vector(-1,0){10}}
\put(70,-15.5){\vector(0,1){10}}
\put(18,0){\makebox(0,0){$\alpha_1$}}
\put(31,0){\makebox(0,0){$\alpha_2$}}
\put(44,0){\makebox(0,0){$\alpha_3$}}
\put(57,0){\makebox(0,0){$\alpha_4$}}
\put(70,0){\makebox(0,0){$\alpha_5$}}
\put(83,0){\makebox(0,0){$\alpha_6$}}
\put(96,0){\makebox(0,0){$\alpha_7$}}
\put(74,-13){\makebox(0,0){$\alpha_8$}}
\put(5,-4){\makebox(0,0){$E_8 \,\ :$}}
\end{picture}
\end{center}

\vspace{18mm}

\noindent
Then, by Lemma \ref{lemma:(3.1)}, one easily sees the following:

\begin{lemma}
\label{lemma:(4.3.1)}
Let $v=v(c_1, \cdots, c_n)$ for $E_n$ $(n=6,7,8)$ with the
orientation of the Dynkin diagram as above. Then 
$\widehat{X}_{\alpha_j}v$ for $j=1, \cdots, n$ are given 
as follows:
\begin{enumerate}
\item[{\rm 1)}] \quad In the case $E_6$,
$$
2\widehat{X}_{\alpha_j}v = \left\{
\begin{array}{lcl}
-ic_jv & & (j=2, 4, 6) \cr
-ic_1v\left(
\begin{array}{rcl}
c_1, \, -c_2 , &\hspace{-2mm}c_3, & \hspace{-2mm}c_4, \, c_5 \cr
&\hspace{-2mm}c_6 &
\end{array}
\right) & & (j=1) \cr
-ic_3v\left(
\begin{array}{rcl}
c_1, \, -c_2 , &\hspace{-2mm}c_3, & \hspace{-2mm}-c_4, \, c_5 \cr
&\hspace{-2mm}-c_6 &
\end{array}
\right) & & (j=3) \cr
-ic_5v\left(
\begin{array}{rcl}
c_1, \, c_2 , &\hspace{-2mm}c_3, & \hspace{-2mm}-c_4, \, c_5 \cr
&\hspace{-2mm}c_6 &
\end{array}
\right) & & (j=5).
\end{array}
\right.
$$
\item[{\rm 2)}] \quad In the case $E_7$,
$$
2\widehat{X}_{\alpha_j}v = \left\{
\begin{array}{ll}
-ic_jv & (j=2, 4, 6, 7) \cr
-ic_1v\left(
\begin{array}{rcl}
c_1, \, -c_2 , &\hspace{-2mm}c_3, & \hspace{-2mm} c_4, 
\, c_5, \, c_6 \cr
&\hspace{-2mm}c_7 &
\end{array}
\right) & (j=1) \cr
-ic_3v\left(
\begin{array}{rcl}
c_1, \, -c_2 , &\hspace{-2mm}c_3, & \hspace{-2mm}-c_4, 
\, c_5, \, c_6 \cr
&\hspace{-2mm}-c_7 &
\end{array}
\right) & (j=3) \cr
-ic_5v\left(
\begin{array}{rcl}
c_1, \, c_2 , &\hspace{-2mm}c_3, & \hspace{-2mm}-c_4, 
\, c_5, \, -c_6 \cr
&\hspace{-2mm}c_7 &
\end{array}
\right) & (j=5).
\end{array}
\right.
$$
\item[{\rm 3)}] \quad In the case $E_8$,
$$
2\widehat{X}_{\alpha_j}v = \left\{
\begin{array}{ll}
-ic_jv & (j=2, 4, 6, 8) \cr
-ic_1v\left(
\begin{array}{rcl}
c_1, \, -c_2 , \, c_3, \, c_4, &\hspace{-2mm} c_5, & 
\hspace{-2mm} c_6, \, c_7 \cr
&\hspace{-2mm} c_8 &
\end{array}
\right) & (j=1) \cr
-ic_3v\left(
\begin{array}{rcl}
c_1, \, -c_2 , \, c_3, \, -c_4, &\hspace{-2mm} c_5, 
& \hspace{-2mm}c_6, \, c_7 \cr
&\hspace{-2mm} c_8 &
\end{array}
\right) & (j=3) \cr
-ic_5v\left(
\begin{array}{rcl}
c_1, \, c_2 , \, c_3, \, -c_4, &\hspace{-2mm} c_5, 
& \hspace{-2mm} -c_6, \, c_7 \cr
&\hspace{-2mm} -c_8 &
\end{array}
\right) & (j=5) \cr
-ic_7v\left(
\begin{array}{rcl}
c_1, \, c_2 , \, c_3, \, c_4, & \hspace{-2mm} c_5, 
& \hspace{-2mm}-c_6, \, c_7 \cr
&\hspace{-2mm}c_8 &
\end{array}
\right) & (j=7).
\end{array}
\right.
$$
\end{enumerate}
\end{lemma}

From this lemma, one sees that, for example in the case $E_7$, 
$c_1$, $c_3$, $c_5$ and $c_4c_6c_7$ are unchanged under the 
action of $\widehat{X}_{\alpha_j}$'s.  
Since $\widehat{X}_{\alpha}$ is the action of the 
field $\widetilde{X}_{\alpha}(z)$ to the $\ccc\{Q/2Q\}$-component, 
one obtains the following:

\begin{thm}
\label{thm:(4.3.1)}
\begin{enumerate}
\item[{\rm 1)}] \quad In the case $E_6$, 
for an arbitrary choice of $c_1, c_3, c_5 \in \{\pm 1\}$, 
the $\ccc$-linear span of
$$ \left\{
v\left(
\begin{array}{rcl}
c_1, \, b_2, & \hspace{-2mm} c_3, & \hspace{-2mm} b_4, \, c_5 \cr
& \hspace{-2mm} b_6&
\end{array} \right)
\,\ ; \,\ 
b_2, b_4, b_6 \in \{\pm 1\} 
\right\}
$$
tensored with $\ccc [ x^{(j)}_r ; 1 \le j \le 6, 
\, r \in \nnn_{\rm odd} ]$ is an irreducible $E^{(2)}_6$-module
of level one.
\item[{\rm 2)}] \quad In the case $E_7$, 
for an arbitrary choice of $c_1, c_3, c_5, c \in \{\pm 1\}$, 
the $\ccc$-linear span of
$$ \left\{
v\left(
\begin{array}{rcl}
c_1, \, b_2, & \hspace{-2mm} c_3, & \hspace{-2mm} b_4, 
\, c_5 , \, b_6 \cr
&\hspace{-2mm} b_7&
\end{array} \right)
\,\ ; \,\ 
\begin{array}{ccl}
{\rm (i)} & &
b_2, b_4, b_6, b_7 \in \{\pm 1\} \cr
{\rm (ii)} & &
b_4b_6b_7 = c
\end{array}
\right\}
$$
tensored with $\ccc [ x^{(j)}_r ; 1 \le j \le 7, 
\, r \in \nnn_{\rm odd} ]$ is an irreducible $E^{(1)}_7$-module
of level one.
\item[{\rm 3)}] \quad In the case $E_8$, 
for an arbitrary choice of $c_1, c_3, c_5, c_7 \in \{\pm 1\}$, 
the $\ccc$-linear span of
$$ \left\{
v\left(
\begin{array}{rcl}
c_1, \, b_2, \, c_3, \, b_4, & \hspace{-2mm} c_5 , & 
\hspace{-2mm} b_6, \, c_7 \cr
&\hspace{-2mm} b_8&
\end{array} \right)
\,\ ; \,\ 
b_2, b_4, b_6, b_8 \in \{\pm 1\} 
\right\}
$$
tensored with $\ccc [ x^{(j)}_r ; 1 \le j \le 8, 
\, r \in \nnn_{\rm odd} ]$ is an irreducible $E^{(1)}_8$-module
of level one.
\end{enumerate}
\end{thm}

The non-irreducibility of the representation space 
$\ccc\{Q/2Q\} \otimes \ccc [ x^{(j)}_r]$ may suggest the 
existence of still bigger symmetry or the action of some
group. 

From the asymptotics of characters of integrable repsentations 
given in \cite{KW1}, one sees that the basic $E^{(1)}_8$-module 
decomposes into the sum of two fundamental representations
of $D^{(1)}_8$.  Actually in the case 3) of the above theorem,
the $\ccc$-linear space of 
$$ \left\{
v\left(
\begin{array}{rcl}
c_1, \, b_2, \, c_3, \, b_4, & \hspace{-2mm} c_5 , & 
\hspace{-2mm} b_6, \, c_7 \cr
&\hspace{-2mm} b_8&
\end{array} \right)
\,\ ; \,\ 
\begin{array}{ccl}
{\rm (i)} & & b_2, b_4, b_6, b_8 \in \{\pm 1\} \cr
{\rm (ii)} & & b_6b_8=c
\end{array}
\right\}
$$
tensored with $\ccc [x^{(j)}_r ] $  is $D^{(1)}_8$-stable and 
$D^{(1)}_8$-irreducible for $c_1, c_3, c_5, c_7, 
c \in \{\pm 1\}$ with respect to the canonical inclusion 
of $D^{(1)}_8$ into $E^{(1)}_8$.

In concluding this note, we remark that the above construction 
gives rise to the product expression of specialized character 
(cf. \cite{Kbook1} \S 10.8) of fundamental representations 
with respect to a particular specialization.
We consider an affine Lie algebra with the simple root system 
$\{\alpha_0, \cdots, \alpha_{\ell}\}$, following 
the enumeration of simple roots from \S 4.8 of \cite{Kbook1}.
Fix a number $s \in \{0 , \cdots, \ell\}$, and consider the
specialization 
$$e^{-\alpha_j} \longmapsto 1 \,\ (j \ne s), \qquad
e^{-\alpha_s} \longmapsto q,$$
which induces an algebra homomorphism
$$F_s \,\ : \,\ \ccc [[e^{-\alpha_0}, \cdots, e^{-\alpha_{\ell}} ]]
\,\ \longrightarrow \,\ \ccc[[ q]] $$
of the associative rings of formal power series.
Then, from Theorems \ref{thm:(4.1.1)}, \ref{thm:(4.2.1)} and 
\ref{thm:(4.3.1)}, we obtain the following expression of 
the specialized character $F_s(e^{-\Lambda}{\rm ch}_{\Lambda})$ 
for a level one dominant integral form $\Lambda$ and 
a specially chosen index $s$:

\begin{cor}
\label{cor:(4.1.1)}
$$
\begin{array}{c|c|c}
\text{affine algebra} & s \, (\text{special index}) 
& \text{specialized character}
\,\ F_s(e^{-\Lambda}{\rm ch}_{\Lambda})
\cr
\noalign{\hrule height0.8pt}
\hfil
A^{(2)}_{2\ell-1} & \ell & 
{\displaystyle 
2^{\ell-1}
\left(\frac{\varphi(q^2)}{\varphi(q)}\right)^{2\ell-1} } \\[3mm]
\hline
A^{(2)}_{2\ell} & \ell & 
{\displaystyle 
2^{\ell}
\left(\frac{\varphi(q^2)}{\varphi(q)}\right)^{2\ell} } \\[3mm]
\hline
D^{(1)}_{\ell} \,\ (\ell : \, {\rm even}) & \ell/2 & 
{\displaystyle 
2^{\ell/2-1}
\left(\frac{\varphi(q^2)}{\varphi(q)}\right)^{\ell} } \\[3mm]
\hline
D^{(2)}_{\ell+1} \,\ (\ell : \, {\rm even}) & \ell/2 & 
{\displaystyle 
2^{\ell/2}
\left(\frac{\varphi(q^2)}{\varphi(q)}\right)^{\ell+1} } \\[3mm]
\hline
E^{(2)}_6  & 4 & 
{\displaystyle 
2^3\left(\frac{\varphi(q^2)}{\varphi(q)}\right)^6 } \\[3mm]
\hline
E^{(1)}_7  & 7 & 
{\displaystyle 
2^3\left(\frac{\varphi(q^2)}{\varphi(q)}\right)^7 } \\[3mm]
\hline
E^{(1)}_8  & 7 & 
{\displaystyle 
2^4\left(\frac{\varphi(q^2)}{\varphi(q)}\right)^8 } \\[3mm]
\noalign{\hrule height0.8pt}
\end{array}
$$
where $\Lambda$ is a dominant integral form of 
level one and $\varphi(q):=\prod\limits^{\infty}_{n=1}(1-q^n)$.
\end{cor}

Actually this corollary holds because ${\rm deg}(e^{\alpha})=0$
(for all $\alpha \in Q/2Q$) and ${\rm deg}(x^{(j)}_r)=r$ with 
respect to this specialization.

We remark that this result is, of course, in coincidence with the 
product expression of characters given in \S 2.3 of \cite{Wbook1}: 
$$
\begin{array}{ccl}
{\rm ch}(\Lambda_0, A^{(2)}_{2\ell}) &=&
{\displaystyle
e^{\Lambda_0}\prod^{\infty}_{j=1}
\prod_{\alpha \in \overline{\Delta}_{\ell, +}}
(1+e^{\pm \frac{\alpha}{2}}e^{-(j-\frac{1}{2})\delta}), 
} \cr
{\rm ch}(\Lambda_0, D^{(2)}_{\ell+1}) &=&
{\displaystyle
e^{\Lambda_0}
\prod^{\infty}_{j=1}(1+e^{-j\delta})
\prod_{\substack{j \in \nnn_{\rm odd} \cr
\alpha \in \overline{\Delta}_{s, +}}}
(1+e^{-j\delta+\alpha}) 
(1+e^{-(j-1)\delta-\alpha}), 
}
\end{array}
$$
where $\overline{\Delta}_{\ell, +}$  
(resp. $\overline{\Delta}_{s, +}$)  is the set of all positive 
long (resp. short) roots of the finite-dimensional Lie algebra 
with the simple root system $\{\alpha_1, \cdots, \alpha_{\ell} \}$.
One easily sees that the specialization of these characters 
gives the same formulas with the above corollary
for $A^{(2)}_{2\ell}$ and $D^{(2)}_{\ell+1}$, 
where $e^{-\delta}=q^2$ in our specialization since 
the coefficient of $\alpha_s$ in the primitive imaginary 
root $\delta$ is equal to $2$.

\end{document}